\documentclass[opre,nonblindrev]{workingpaperJJB}

\OneAndAHalfSpacedXI 

\usepackage{endnotes}
\let\footnote=\endnote
\let\enotesize=\normalsize
\def\notesname{Endnotes}%
\def\makeenmark{\hbox to1.275em{\theenmark.\enskip\hss}}
\def\enoteformat{\rightskip0pt\leftskip0pt\parindent=1.275em
  \leavevmode\llap{\makeenmark}}

\usepackage{algorithm}
\usepackage{algpseudocode}
\usepackage[section]{placeins}
\usepackage{multirow}
\usepackage{diagbox}
\usepackage{array}
\usepackage[justification=centering]{caption}
\usepackage[justification=centering]{subcaption}
\usepackage{bbold}
\usepackage{appendix}
\usepackage{xcolor}
\usepackage{natbib}
\setcitestyle{round}

\TheoremsNumberedThrough     
\ECRepeatTheorems

\EquationsNumberedThrough    

\begin{document}


\RUNAUTHOR{Mintz and He}

\RUNTITLE{Model Based Reinforcement Learning for Personalized Heparin Dosing}

\TITLE{Model Based Reinforcement Learning for Personalized Heparin Dosing}

\ARTICLEAUTHORS{%
\AUTHOR{Qinyang He and Yonatan Mintz}
\AFF{Department of Industrial and Systems Engineering, University of Wisconsin -- Madison, Madison, WI 30332, \EMAIL{\{qhe57,ymintz\}@wisc.edu}}
} 

\ABSTRACT{%
One of the key challenges in sequential decision making is optimizing systems safely in the case of partial information. While much of the existing work has focused on addressing this challenge in the case of either partially known states or partially known system dynamics, it is further exacerbated in cases where both states and dynamics are partially known. For instance  the setting of computing heparin doses for patients fits this paradigm since the concentration of heparin in the patient cannot be measured directly and the rates at which patients metabolize it vary greatly between individuals. While many approaches proposed to resolve the challenge in this setting are model free, they require complex models that are not transparent to decision makers, and are difficult to analyze and guarantee safety. However, if some of the structure of the dynamics is known, a model based approach can be leveraged to provide safe policies with practical empirical performance and theoretical worst case guarantees. In this paper we propose a model based framework to address the challenge of partially observed states and dynamics in the context of designing personalized doses of heparin. We use a predictive model based on pharmacokinetics (the study of how the body effects substances through absorption, distribution, and metabolism) parameterized individually by patient, and infer the current concentration of heparin and predict future therapeutic effects taking into account different patients' characteristics. We formulate the patient parameter estimation problem in to a mixed integer linear program and show that our estimates are statistically consistent. We leverage this model by developing an adaptive dosing algorithm that outputs asymptotically optimal dose sequences based on a scenario generation approach, this approach is also capable of ensuring that the required heparin doses are maintained within a safe level. We validate our models with numerical experiments by first comparing the predictive capabilities of our model against existing machine learning techniques and demonstrating how our dosing algorithm can keep patients' related medical tests within a therapeutic range in a simulated ICU environment. Our results show that our methods are capable of maintaining patients in therapeutic range for 87.7\% of the treatment time  as opposed to existing weight based protocols that can only do so for 55.6\% of the treatment time. 

}%


\KEYWORDS{personalized healthcare, sequential decision making, integer programming, predictive analytics}

\maketitle
%
\section{Introduction}


One of the biggest challenges facing decision makers when making repeated decisions under uncertainty is when their actions are time sensitive and costly. In the operations literature many models proposed for these challenges have assumed that either the only uncertainty in the problem is in the state of the system (i.e. it is partially observed) but the dynamics of the system are known \citep{unknown_state,POMDP_1}, or there is uncertainty about the dynamics of the system but the state is fully observed \citep{unknown_dynamics,sutton2018reinforcement}. However, in many scenarios uncertainty exists both in the dynamics and states of the system.  For instance, in the setting of Unfractionated Heparin dosing. Unfractionated Heparin (often times simply referred to as heparin) is a common anticoagulant used in intensive care units (ICUs) \citep{brunet2008pharmacodynamics}. In this setting, a clinician must compute the most effective sequence of heparin doses to administer to a patient over the course of the treatment while ensuring they stay in therapeutic levels of blood coagulation.  However, the clinician can never directly observe the concentration of heparin in the patient's body and can only use noisy lab measurements to adjust their dosing sequence. Moreover, the half-life of heparin in a patient's body (that is the amount of time it would take a patient to metabolize through half of a dose) can vary from 60-90 minutes in most cases, but depending on the size of the dose and patient characteristics can be as short as 30 minutes and up to 2.5 hours \citep{cook2010anticoagulation,heparin_problem}. In other words, the clinician only has partial information on the health state of the patient that can only be learned with noisy observations and does not fully know the dynamics of how quickly the drug will be metabolized in the patient's body.


One set of approaches that have been proposed to address this level of uncertainty in heparin dosing are model-free approaches \citep{deeparch,Nemati2016OptimalMD}. These methods treat lab test results and chart data as a fully observed state and approximate the dynamics of the problem using complex models such as neural networks that can capture a large variety of functional relations. In principle, model free approaches should adapt to any dynamic structure given enough data and sufficient exploration and have been used to great effect in non-healthcare applications such as gaming and robotics \citep{li2017deep}. However, these approaches may not be suitable for time sensitive settings such as heparin dosing as in practice they take a long time to converge, are sensitive to missing and noisy data \citep{model_free}, and are not generally interpretable to decision makers. Moreover, it is difficult to ensure that the underlying system being controlled by these model-free approaches is staying within a designated safety range. However, in many contexts the structure of the system could be well known from previous research, and although the true state and dynamics may not be known this structure can be exploited to design safe and effective algorithms. In particular, for a drug such as heparin that has well studied pharmacology, the structure of existing pharmacological models can be leveraged to provide effective dosing protocols. 

In this paper, we propose a model based framework that is capable of optimizing costly decisions in the context of partially known system dynamics and  partially observed system states in a safe manner. We focus our modeling and analysis on the case of personalized heparin dosing, though portions of our framework could be applied in more general settings. Our framework is built from two key components (i) an individual level model of each patient receiving treatment and (ii) an optimization problem that computes the appropriate dose sequence for each patient. The patient model describes how each individual is metabolizing the heparin in their body, their level of heparin concentration, and how their blood coagulation reacts to the level of heparin. Using patient level lab and chart data, we can use our model to compute several possible scenarios that capture the uncertainty in the patient's dynamics and health state, and evaluate the likelihood of these scenarios occurring. Using these scenarios our dose sequence optimization will compute a sequence of doses that ensures patients will remain within a therapeutic range and avoid dangerous levels of blood coagulation, while accounting for uncertainty in the model parameters. In practice, this framework will be implemented itteratively and adaptively, as more data is collected from the patient the likelihood estimates of the scenarios will be recalculated and so will the dosing sequence.


\subsection{Unfractionated Heparin Dosing}
\label{sec:ufh_dose}
Administering  heparin is a common practice in ICUs because it prevents blood clotting in patients undergoing hemodialysis \citep{shen2012use}. In general, when a patient is ordered a heparin treatment they will be first given an initial dose called a bolus dose, and then every four to six hours the medical team will administer a new dose at varying levels in response to the patient's lab tests \citep{hirsh2001guide}. The most common lab measure used is activated partial thromboplastin time (aPTT or just PTT) which is a measure of blood coagulation \citep{shen2012use,hirsh2001guide}. Generally, physicians modify their doses to target a therapeutic aPTT level that is 1.5 times higher then the patient's baseline level, or alternatively ensure that patient aPTT measures remain between 1.5 to 2.5 times their baseline level \citep{hirsh2001guide}.

However, since heparin is fast acting, improperly dosing patients with heparin introduces several serious risks. Specifically, overdosing heparin can lead a patient to hemorrhage while underdosing can result in clotting \citep{landefeld1987identification}. Although heparin has been the most commonly used anticoagulant in the U.S.,  one-third of patients that are given heparin treatment are misdosed \citep{secretariat2009point}. This is because each individual metabolizes heparin at different rates and inferring each individual patient's metabolism rates based on their demographics (such as age and gender) or past medical records is very challenging. To address this challenge we design a precision medicine framework that addresses the heparin dosing problem by personalizing each patient's dosing policy based on their individual data.

To date, there have been several proposed data driven heparin dosing policies in the medical literature. One guiding principle is to apply an hourly dose within specific range. The range proposed for the initial bolus is 2,000-4,000 international standard units (IU) of heparin \citep{davenport2011optimization} and the range for the hourly rate is from 500-2,000 IU/h or more \citep{wilhelmsson1984heparin}. However, these ranges are not very reliable since they do not adjust for patient specific features and do not eliminate misdosing \citep{shen2012use}. Other common protocols for heparin dosing are weight-based \citep{weight1,weight2}, where the initial bolus and hourly rate are decided using the patient's weight or body mass index (BMI). However, weight based policies are not significantly better than fixed dose policies in maintaining  patients in the targeted aPTT range \citep{shen2012use}.

With the advent of Electronic Health Records (EHR) and large amounts of patient level data,  modern data-driven approaches have been proposed to address the challenges of heparin dosing. These methods include supervised learning methods such as multivariate logistic regression \citep{multilogit} that aims to predict whether the patient is at sub-therapeutic or super-therapeutic level (their aPTT measure is below or above the safe therapeutic range respectively) given an initial dose. Additionally reinforcement learning (RL) methods \citep{Nemati2016OptimalMD,deeparch}  have also been propose. In general these methods are designed to learn a single dosing policy from dosing trial data in EHR. In particular \citet{deeparch} develop a complex deep learning architecture that combines deep belief network to encode system states and a softmax regression output layer for dose prediction. All these methods are trained on full patient records and assume all patients follow the same heparin metabolism dynamics. Thus they can be thought of as providing a one size fits all dosing policy for patients. Therefore, even though the overall empirical performance of these methods seems strong, the resulting policy might not be suitable for a specific patient. In this paper, we will consider a personalized heparin dosing setting where individual patient data will be used to guide dosing decisions.

\subsection{Related Literature}
The methods we develop in this paper are related to several streams of research within the operations community that include safe reinforcement learning, partially observed Markov decision processes (POMDPs), personalized dosing, and sequential decision making.

Our work in this paper contributes to the larger stream of safe reinforcement learning that has been a key area of interest in the literature. In this setting a decision maker must solve a sequential decision making problem with only partial information on the dynamics of the model, while ensuring that the state of the system does not enter an ``unsafe'' subset of the state space. One approach that has been proposed to address this problem has been to include an explicit safety factor in the stage costs of the sequential decision making problem, that is a term that estimates the probability the state will transition into an unsafe state given the current control action \citep{saferl_1,saferl_2, saferl_3}. In addition to modifications of the cost function, a different stream of literature has considered explicitly modifying the exploration process of states to mitigate the risks of reaching an unsafe state due to random exploration \citep{saferl_6, saferl_4, saferl_5}. Concepts of safety have also been applied in bandit problems \citep{safebandit1, safebandit2}, here the notion of risk comes from a set of constraints set on the rewards of the arms and less from a predefined unsafe set. The main challenge of applying these methods in the case of heparin dosing is that in order to achieve strong performance they still require a large amount of exploration, that could translate into delayed effective treatment and adverse effects for patients. Our proposed framework builds from the exploration and constrained based approaches that can use patient data more efficiently, and in our experiments converges quickly to effective dosing recommendations.

Another class of models related to our modeling framework is that of Partially Observable Markov Decision Processes (POMDPs). In this class of problems, the decision maker is solving a Markov Decision Process (MDP) with only partial observations of the true system state. The canonical methods for optimizing POMDPs involve converting the POMDP into a belief MDP \citep{POMDP_1}, that is a full information MDP where the new states of the system encode the probability the true system state is some value, and then solve the belief MDP using dynamic programming \citep{POMDP_2}. While this method is computationally tractable for small POMDPs, the resulting belief MDP is often too large to optimize effectively \citep{POMDP_3}. Thus in many large scale setting  approximate dynamic programming techniques have been proposed to approximate optimal policies \citep{POMDP_4,jiang2015approximate,wang2022optimal}. 
Our setting differs from the POMDP setting in that the transition function is also partially observable since each patient has varying metabolism rate that makes the belief state conversion challenging to implement.

Both reinforcement learning \citep{yu2021reinforcement} and POMDPs have been applied in the setting of personalized healthcare \citep{keskinocak2020review}. Several reinforcement learning methods have been proposed for personalized medicine including deep reinforcement learning \citep{pdosing_2,pdosing_1}, off policy learning \citep{wang2022reliable}, and some model based methods \citep{lee2015applying,skandari2021patient}. However,developing reinforcement learning solutions to increase the safety and robustness of learned strategies in healthcare remains a key challenge \citep{pdosing_3}.
%
%
POMDPs have also been used to study personalized dosing in applications such as heart disease \citep{pdosing_4}, stroke \citep{pdosing_6} and Parkinson's disease \citep{pdosing_5}. POMDP models have also been used in population health and disease screening and disease management \citep{pdosing_7,bonifonte2022analytics,wu2022optimizing,bertsimas2020personalized,zhang2012optimization,hajjar2023personalized,zhang2012optimization}. \cite{shi2021timing} incorporated a personalized readmission prediction model to manage patient flow in a hospital. The methods we propose in this paper expand on these methods to a healthcare setting where both states and dynamics are partially observed.

Another stream of literature closely related to the framework proposed in this paper use explicit dynamic models. \citet{mintz2017behavioral} proposed a behavioral analytics framework to optimized incentives given by a single coordinator to a multi-agent system. This sequential decision making framework has some similarities with the heparin dosing problem but from an incentive optimization perspective. However, the heparin dosing problem has more pronounced safety concerns then mobile weight loss interventions. \cite{dogan2021regret} developed learning-based policies for model predictive control that require parameter estimation for the dynamics. \citet{lee2018outcome} develop a predictive dose-effect model for diabetes management that establishes a direct relationship between drug dose and drug effect and also a treatment planning optimization model. Our framework builds on these methods, and extends them to the case of personalized heparin dosing.

\subsection{Contributions}
In this paper we develop a model based framework for personalized heparin dosing under partially observed states and dynamics. By developing this framework, we provide four major contributions:
\begin{enumerate}
    \item We propose a patient level model for how heparin is metabolized by each individual that can be incorporated into an optimization framework. To the best of our knowledge we propose one of the first piece-wise linear treatment-effect model for heparin dynamics based on the Michaelis-Menten equations \citep{MM1}. This model captures the  individual rate at which patients metabolize heparin and is parameterized by patient specific parameters. In contrast to machine learning approaches that train a single model for different patients, our model allows us to personalize patient state prediction and future dose planning in a more transparent and interpretable way by expressing the relationship between heparin level and doses applied.
    
    \item  We propose an estimation technique that can simultaneously estimate both unknown individual states and dynamic parameters of each patient. Unlike most dynamic models used in the operations literature where there is either unknown state or unknown transition function, our work addresses uncertainty in both patients' state and state transition dynamics. Through a joint maximum likelihood estimation approach, we formulate an mixed integer optimization problem  (MIP) that solves the patient state estimation and parameter learning problem simultaneously. In addition, we propose a decomposition scheme that can quickly solve the training problem. We also show that parameters estimated using this method are statistically consistent.

    \item We develop a dynamic heparin dosing algorithm that can provide individual level dosing sequences that maintain a patient's aPTT within a safe range. Instead of training a single model for all the patients, our approach models the interaction between the clinician's medical decision and each individual patient's reaction in the context of heparin treatment. Our adaptive dosing algorithm adjusts each patient's future doses based on their individual aPTT measurement history. We also account for safety in our algorithm design by estimating the uncertainty in both state and dynamic parameter estimation. We provide a statistical guarantee that our algorithm is asymptotically optimal, which means that as more patient observations are collected the algorithm will provided optimal doses in line with the patient's true condition.

    \item We conduct two sets of numerical experiments to validate our modeling and dosing algorithm using the MIMIC III \citep{mimic3} data set. The first set of experiment compares the predictive accuracy of our model against other machine learning approaches and the results show that our dynamic model outperform other state-of-the-art methods in predicting future patient aPTT levels. The second set of experiment tests how our adaptive algorithm performs in maintaining patient aPTT in a safe range. Based on a simulation with 25 patients, regardless of their initial conditions, the best implementation of our algorithm on average keep patients in the safe therapeutic range 87.7\% of time over the course of their treatment as opposed to the existing weight based protocol that does so for 55.6\% of the treatment time.

\end{enumerate}

\subsection{Paper Overview}
In Section \ref{sec:model_desc}, we describe our model for heparin dynamics. We introduce the nonlinear pharmacokinetic model in medical literature that characterizes the rate at which a patient metabolizes heparin and implement a piece-wise linear approximation of this model so that it can be used in commercial optimization software while preserving many of its unique properties. We then formulate the dynamics of aPTT and its interaction with the heparin concentration. We combine these formulations into a single heparin dynamic model parameterized by patient-specific characteristics.

In Section \ref{sec:m_param_est}, we develop methods for estimating the patient-specific parameters of the model using available aPTT observations through a joint parameter maximum likelihood estimation approach. We show that the estimators obtained by this approach are consistent. We reformulate the parameter estimation problem into a MIP. To tackle the computational complexity of this MIP formulation, we  propose a generalized Benders Decomposition algorithm to learn the unknown dynamics and states with a faster computational time.

In Section \ref{sec:dose_opt}, we explore methods to design personalized optimal dosing policy for each patient. We formally define the dose optimization problem and propose an adaptive dosing algorithm based on a scenario generation approach. We then prove the asymptotic optimality of our algorithm.

We conclude with Section \ref{experiments} where we conduct two sets of numerical experiments to test our framework using the MIMIC III data set. First, we compare the predictive accuracy of our dynamic model to other existing methods. Second, we evaluate the effectiveness of our adaptive dosing algorithm against existing weight based protocols using a simulation study.

\section{Model Description}
\label{sec:model_desc}
In this section we develop the dynamic model that describes the trajectory of heparin in the patient body and its relationship with aPTT measurements. For our model we define  $x_t \in \mathcal{X}, u_t \in \mathcal{U}, y_t\in \mathcal{Y} $ as the concentration of heparin, the heparin dose, and the aPTT of the patient at hour $t = 0,...,T$ of the patient's stay in the ICU. We assume $\mathcal{X}, \mathcal{U}, \mathcal{Y} \subset \mathbb{R}_+$ are closed intervals. We also assume that the initial level of heparin in the patient's body is zero that is $x_0 = 0$. This implies no heparin has been administered prior to their admission into the ICU. On the other hand we assume that the initial aPTT is unknown but is close to its baseline value. Our goal is to model the dynamics of heparin as a parametric grey box dynamic model of the form:
\begin{equation}
   \begin{aligned}
   y_t &= f(x_t,y_{t-1},\theta_y),\\
   x_t & = g(x_{t-1},u_t,\theta_x).
   \end{aligned}
\end{equation}

Here $f: \mathcal{X}\times \mathcal{Y} \times \Theta_y \rightarrow \mathcal{Y}$ and $g: \mathcal{X} \times \mathcal{U} \times \Theta_x \rightarrow \mathcal{X}$ are dynamics functions of known form and $\theta_y,\theta_x$ are unknown parameters that characterize the difference in dynamics between patients. We will generally assume  the parameter sets $\Theta_x,\Theta_y$ are compact, however through our analysis and modeling we will enumerate additional assumptions about their structure to ensure model identifiability. It is important to note that since we assume the form of $f,g$ are known there are some function parameters the we assume are known \textit{a priori} as opposed to the initially unknown parameters that are all captured by $\theta_y,\theta_x$. Specifically, let $\theta_x = [\alpha, k]$  correspond to the rates at which heparin is metabolized, and $\theta_y = [b, y_0, y_{b0}, y_b]$  correspond to the coefficient of heparin influencing aPTT, the initial aPTT level, the temporary elevated baseline aPTT, and the long run homeostasis aPTT of the patient. The detailed meaning of these parameters will be explained later in this section when describing the functional form of the dynamics and their relationship with existing medical models. These parameters are considered to be unique to each individual patient and thus allow us to personalize patient state prediction and heparin dosage. In other words, two different patients in the ICU will have different values of $\theta_x$ or $\theta_y$ depending on their biological or other health related factors that will cause their treatment plans to differ. 

\subsection{Medical Models of Heparin Metabolism}

Several models have been proposed in the medical literature for how drugs and toxins metabolize in the body \citep{meta1, meta2}. In the terminology of pharmacokinetics (the study of how the body effects substances through absorption, distribution, and metabolism) a key component of these models is the elimination rate,  the rate at which the drug is removed from the body through metabolism and absorption \citep{jambhekar2009basic}. Heparin is particularly hard to model because its elimination rate is proportional to its concentration. In the medical literature this property is known as a nonlinear elimination rate.  Common models describing heparin dynamics are continuous in time making them challenging to incorporate into an optimization framework.

To understand these models consider the following modeling example. Suppose a patient is admitted to the ICU at time 0 and a physician determines that they need to be treated using heparin. Using our previous notation, let $x_t$ be the concentration of heparin at time $t$ in the patient's body and $x_1$ be the concentration of the drug after the administration of the initial bolus dose. 
%
Two common elimination models used in pharmacology are zero-order elimination and first-order elimination \citep{jambhekar2009basic}. Zero-order elimination implies that the concentration of the drug decreases linearly over time that is $\frac{dx_t}{dt} \propto 1$. Alternatively first-order elimination implies exponential decay of the concentration of the drug, hence it can be described by $\frac{dx_t}{dt} \propto x_t$.

One of the challenging aspects of appropriately dosing heparin is that its elimination is nonlinear, that is they do not conform to either first-order or zero-order elimination rates. In fact the rate of elimination of heparin varies from acting more like first-order to zero-order depending on the concentration of heparin. A common model proposed to characterize these kinds of nonlinear pharmacokinetics is known as the Michaelis-Menten equations (MM) \citep{MM1}. Originally proposed to model the production of a product from some kind of enzymatic reaction, these equations essentially model a faster rate of elimination with higher concentration. Existing literature has shown that heparin elimination dynamics are modeled extremely well by the MM dynamics \citep{MM2}  and hence provide a good basis for a grey box model. The MM dynamics are characterized by the following differential equation:
\begin{equation} \label{eq:mm_dynam}
    \frac{dx_t}{dt} = V_{max}\frac{x_t}{x_t + K}.
\end{equation}
Here $V_{max}$ stands for the maximum rate of elimination and $K$ represents the value at which the maximum rate is halved. The intuition for these dynamics comes from their limiting behavior with respect to $x_t$. Note that when $x_t << K$ the derivative becomes approximately proportional to $x_t$ implying an exponential rate of decay, i.e. first-order elimination rates. However, when $x_t \geq K$ the rate simplifies to a rate proportional to $V_{max}$ that does not depend on $x_t$ implying zero-order elimination. While the concentration is large the rate of elimination is a fast linear rate, but once the concentration decreases it starts decaying exponentially towards a level of zero. 

Note that the MM equation implies dynamics that are continuous in time. While this may be accurate to describe the process of drug elimination, this is difficult to use in decision making settings. In particular, since clinicians monitor the administration of heparin periodically, they are operating in a discrete time regime and not continuous time. Therefore to create a method that can provide relevant insights, we need to discretize the dynamics with respect to time. To do this requires using the explicit solution to the MM equations derived by \cite{MM3}:
\begin{equation} \label{eq:mm_solution}
    x(t) = K ~ W\Big(\frac{x_1}{K} \exp\big(\frac{x_1}{K} - \frac{V_{max}}{K}t \big) \Big).
\end{equation}

Where $K,V_{max}$ are as previously defined in Equation \eqref{eq:mm_dynam}, $x_0$ represents the initial dose, and $W: \mathbb{R} \rightarrow \mathbb{R}$ is the Lambert $W$ function (that is $W(x\exp(x)) = x$). 

\subsection{Piece-wise Linear Approximation of Heparin Dynamics} 
While the MM equations are useful for describing non-linear elimination they are difficult to use in optimization frameworks. This is because as shown in \eqref{eq:mm_solution} the closed form solution of the equation is both non-trivial and transcendental. This means identifying the relevant model parameters and optimizing dosing policies cannot be done efficiently using commercial optimization software. Existing methods for identifying these parameters have been developed \citep{heparin2} however, these do not guarantee strong statistical properties such as consistency and the resulting model is still difficult to incorporate into a discrete time decision making framework. Hence usually these dynamics have been fully discretized for optimizing dosing policies \citep{Nemati2016OptimalMD} and estimated form data using enumeration or complex least squares methods \citep{heparin_1}.

Instead of full discretization, we can preserve many of the useful qualities of the MM dynamics by implementing a piece-wise linear approximation of the dynamic equations for discrete time. To see how this can be achieved consider the following dynamics equation:

\begin{equation}
\label{eq:linearized_mm}
    x_{t+1} = h(x_t; \alpha, k) = \begin{cases}\alpha x_t, \quad \text{for } x_t \leq \frac{k}{1-\alpha}, \\ x_t-k, \quad \text{for } x_t > \frac{k}{1-\alpha} \end{cases}.
\end{equation}

\begin{figure}[h]
    \centering
    \includegraphics[width=0.5\textwidth]{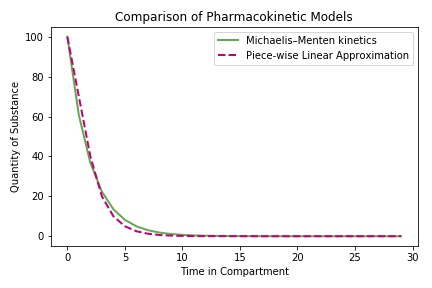}
    \caption{Comparison of MM trajectory (in solid line) and piecewise linear approximation (in dashed line). The x-axis corresponds to time in that patient's body (also called compartment) and the y-axis corresponds to the quantity of the substance.}
    \label{fig:kinetics_comapre}
\end{figure}

Here the two linear components correspond to the limiting behaviors of the MM equation. The parameter $\alpha$ corresponds to the first-order elimination rate of the substance while the parameter $k$ corresponds to the zero-order elimination rate. Observe that when the concentration is high the substance undergoes first-order elimination until reaching a critical value below which it goes through exponential decay, hence the elimination rate is dependent on the concentration. The choice of the threshold at $\frac{k}{1-\alpha}$ is to ensure that the dynamics are continuous. This mimics the limiting behaviors of the MM equation that are fundamental for its motivation. However, using this formulation results in an abrupt change between first-order and zero-order elimination instead of a smooth transition between them. That said, this formulation can still provide close approximations to the MM trajectory. Consider for instance the example in Figure \ref{fig:kinetics_comapre}, we plot two trajectories one generated by MM and one by the approximation given in \eqref{eq:linearized_mm}. Note that as desired the two approximations essentially match in the limiting cases (close to times 0 and as $t$ gets large), with the greatest gap being where the threshold of the linearization is; however, this gap is still quite small. This motivates that this formulation should be able to closely capture the dynamics of heparin elimination. To complete our model of heparin concentration dynamics, we need to acount for the dose administered by the clinician at time $t$, which using our notation is $u_t$. Since the dose adds linearly to the amount of heparin in the body the complete dynamics of heparin are given by $x_{t+1} = h(x_t;\alpha,k) + u_t$, where $h$ is as defined in \eqref{eq:linearized_mm}.

\subsection{Interaction with Lab-Measurements}
In addition to modeling the dynamics of heparin elimination (also referred to as the kinetics of heparin), since we do not obtain direct measurements of the concentration of heparin in the patient's body we need to use aPTT measurements \citep{heparin_aPTT} to estimate it. This requires modeling the interaction between these observed measurements and the concentration of heparin. As observed by \cite{aPTT_linear} the relationship between aPTT and heparin is approximately linear, thus we can consider a set of linear dynamics to describe the trajectory of a patient's aPTT. In particular we propose the following set of dynamics:
\begin{align}
y_{t+1} &= \gamma_1 (y_t - y_{bt}) + y_{bt} + bx_{t+1}, \\
y_{b,t+1} &= \gamma_2 y_{bt} + \gamma_3 y_{t} + \gamma_4 y_b.
\end{align}

The structure of these dynamics is such that over the long run the aPTT of the patient tends to a homeostasis level $y_b$. Here $y_t$ represents the aPTT, $y_{bt}$ represents the medium time length stable state aPTT level, $y_b$ represents the overall homeostasis aPTT of the patient, and $x_t$ as before represents the concentration of heparin. The parameters $\gamma_1,\gamma_2,\gamma_3,\gamma_4 \in (0,1)$ represent the rate at which values decay to their base level. Essentially $\gamma_1$ indicates how current aPTT level returns to its medium term baseline, and  insuring $\gamma_2 + \gamma_3 + \gamma_4 = 1$ shows that $y_{bt}$ is elevated as a geometric mean of previous aPTT levels and the baseline level. Through validation we found that the values of these parameters do not vary to widely between individuals and can be treated as constants. The parameter $b>0$ reflects how the concentration of heparin impacts the aPTT, since this parameter does vary between patients we treat it as initially unknown and it must be estimated from data.

\section{Model Parameter Estimation}
\label{sec:m_param_est}
In this section we consider estimating the unknown parameters in the model described in Section \ref{sec:model_desc} using available data. As mentioned in Section \ref{sec:ufh_dose} heparin is dosed in the following way in practice: at the beginning of the treatment, the clinician will apply a large initial bolus dose to raise the patient's aPTT to therapeutic levels, and then provide an upkeep dose to make sure this level remains within the therapeutic range. This upkeep dose is then updated as the clinician makes periodic observations of the aPTT through lab tests every 4-6 hours. However, due to the nature of these lab tests, there may be hours where the aPTT measurements are not properly recorded causing missing data in the patient's record in addition to significant measurement noise. Therefore, our proposed estimation method will need to address both missing observations and account for the influence of observation errors on estimates of the parameters.

One approach proposed for estimating parameters and missing measurements is a joint parameter maximum likelihood estimation (MLE) approach \citep{embretson2013item}. Suppose the patient has been in the ICU for $T$ hours and we would like to use aPTT observations from these hours to estimate the parameters of the model. To use this kind of approach we need to assume a particular sample noise model. Let $\tilde{y}_t$ be the aPTT observed by the clinician at time $t$. We will assume that these observations relate to the underlying model through $\tilde{y}_t = y_t + \epsilon_t$ where $\epsilon_t$ are i.i.d. random variables such such that $\mathbb{E}\epsilon_t = 0$ and $\mathbb{E}\epsilon_t^2 < \infty$. Let $\{u_t\}_{t=0}^T$ denote previous heparin doses up to time $T$. Using this noise model we can write our MLE problem as the following:

\begin{equation}
\label{eq:mle_prob}
   (\hat{\theta}_{x,MLE}, \hat{\theta}_{y,MLE}) = \argmax_{\theta_x,\theta_y \in \Theta_x \times \Theta_y} p\big(\{\tilde{y}_t\}_{t\in n(T)} | \theta_x,\theta_y, \{u_t\}_{t=0}^T \big).
\end{equation}

Where $n(T) \subset \{1,...,T\}$ represents the set of time periods where an aPTT observation was collected, $p:\Theta_x,\Theta_y \rightarrow \mathbb{R}$ represents the joint p.d.f., and $\hat{\theta}_{x,MLE}, \hat{\theta}_{y,MLE}$ are the MLE estimates of the unknown parameters.  

In this section we show that under reasonable modeling assumptions on the noise distribution, the optimization problem in Equation \eqref{eq:mle_prob} can be formulated as a MIP. We also discuss how this MIP can be solved efficiently using a decomposition scheme that can be implemented with commercial solvers. Finally, We prove the statistical properties of these estimates in particular we show that $\hat{\theta}_{x,MLE}, \hat{\theta}_{y,MLE}$ are consistent estimates of the ground truth parameters in a Bayesian sense. 

\subsection{MILP Formulation of MLE}
\label{sec:milp_form_mle}
To formulate \eqref{eq:mle_prob} as a MILP we first need to rewrite the MLE problem in terms of the model from Section \ref{sec:model_desc}. Using this model we can expand the likelihood as follows:

\begin{equation}
    \begin{aligned}
    &p\big(\{\tilde{y}_t\}_{t\in n(T)} | \theta_x,\theta_y \big) = p(\{\tilde{y}_t\}_{t\in n(T)} | \{y_t\}_{t=0}^T, \{y_{bt}\}_{t=0}^T,\{x_t\}_{t=0}^T,\theta_x,\theta_y,\{u_t\}_{t=0}^n), \\
& = \prod_{t\in n(T)}p(\tilde{y}_t|  \{y_\tau\}_{\tau=0}^t,\{y_{b\tau}\}_{\tau=0}^T,\{x_\tau\}_{\tau=0}^t,\theta_x,\theta_y,\{u_\tau\}_{\tau=0}^t), \\
&= \prod_{t\in n(T)}p(\tilde{y}_t|  y_t) \prod_{t = 1}^T p(y_t | y_{t-1},y_{bt},x_t,\theta_y)p(x_t|x_{t-1},u_t,\theta_x).
\end{aligned}
\end{equation}

Using standard techniques, we can take the log of this expression to write likelihood in an additive from as follows:
\begin{equation}
    \label{eq:add_lik}
    \sum_{t\in n(T)} \log p(\tilde{y}_t|y_t) + \sum_{t=1}^T \log p(y_t|y_{t-1},y_{bt}, x_t, \theta_y) + \sum_{t=1}^T\log p(x_t| x_{t-1}, u_t, \theta_x).
\end{equation}

Note that using the model from Section \ref{sec:model_desc}, the conditional density functions $p(y_t|y_{t-1},y_{bt}, x_t, \theta_y), p(x_t| x_{t-1}, u_t, \theta_x)$ are degenerate and can be represented as a set of constraints. So we can express \eqref{eq:mle_prob} as the following constrained optimization problem:
\begin{subequations}
\label{eq:constraind_mle}
    \begin{align}
    \max\limits_{\{x_t,y_t,y_{bt}\}_{t=0}^T,y_b,\alpha,k,b} & \sum_{t\in n(T)} \log p(\tilde{y}_t|y_t)\\
\text{subject to } &\\
&  x_{t+1} = u_{t+1} + \begin{cases}\alpha x_t \quad \text{for } x_t \leq \frac{k}{1-\alpha} \\ x_t-k \quad \text{for } x_t > \frac{k}{1-\alpha} \end{cases} & \forall t \in \{0,...,T-1\}, \label{eq:hep_dynam_const} \\
& y_{t+1} = \gamma_1 (y_t - y_{bt}) + y_{bt} + bx_{t+1} & \forall t \in \{0,...,T-1\}, \label{eq:ptt_dynam_const}\\
&y_{b,t+1} = \gamma_2y_{bt} + \gamma_3 y_{t} + \gamma_4 y_b & \forall t \in \{0,...,T-1\},\label{eq:ptt_baseline_const}\\
&x_0 = 0, x_t \in \mathcal{X}, y_t,y_{bt},y_b \in \mathcal{Y} & \forall t \in \{0,...,T\}, \\ 
& \alpha \in \mathcal{A}, k\in \mathcal{K},b \in \mathcal{B}.
    \end{align}
\end{subequations}

Where $\mathcal{K},\mathcal{B} \subset (0,\infty)$ are compact subsets for the possible values of $k$ and $b$. As it stands, this formulation is currently not in a form that can be used in commercial solver. This is because first we have the piece-wise linear dynamics of $x_t$ and second the bi-linear terms introduced by $\alpha x_t$ and $b x_{t+1}$. For this analysis we will need to make the following technical assumptions about the model and its parameters.
\begin{assumption}
\label{as:alpha_set_assump}
For first order decay rate $\alpha \in \mathcal{A}$, the set $\mathcal{A}\subset (0,1)$ is finite, that is $|\mathcal{A}| = m < \infty$. Moreover there exists $\epsilon_\alpha$ such that $\forall \alpha \in \mathcal{A}, \alpha > \epsilon_\alpha$. Moreover, the intervals $\mathcal{X}$ and $\mathcal{Y}$ are non-empty and $0 \in \mathcal{X}$.
\end{assumption}

 In practice, the first portion of this assumption generally holds since half lives of substances are only measured up to a certain accuracy for dosing decision purposes (usually to the 15-30min mark). Moreover, generally the minimum possible half life will not be instantaneous. and this assumption will enable us to reformulate the bi-linear terms using a product of binary and continuous variables. The second part of the assumption ensures \eqref{eq:constraind_mle} has a feasible solution.

\begin{assumption}
Let $p_\epsilon(\cdot)$ be the p.d.f. of the distribution of $\epsilon_t$, then $\log p_\epsilon(\cdot)$ is concave and can be described by mixed integer linear constraints and objective terms. 
\end{assumption}

The assumption on the concavity of the log p.d.f. is common in statistical estimation problems \citep{boyd2004convex} and the second half of the assumption ensures that the log likelihood can be used as the objective function of a MILP. This assumption is not restrictive since many common distributions such as the Laplace and exponential distributions follow this assumption in addition to piecewise linear concave distributions.

With the above two assumptions, we are now able to reformulated the problem. We first handle the bilinear terms in the constraints. Consider the coefficient $b$.

\begin{proposition}
\label{prop:z_bx}
Let $z_t = b\cdot x_t$ and $c = b\cdot k$, then Constraints \eqref{eq:hep_dynam_const}--\eqref{eq:ptt_baseline_const} can be reformulated as:
\begin{equation}
\label{eq:reform_bilinear}
    \begin{aligned}
    &z_{t+1} = bu_t + g(z_t) = b u_t + \begin{cases}\alpha z_t, \quad z_t \leq \frac{c}{1-\alpha} \\ z_t - c, \quad z_t > \frac{c}{1-\alpha} \end{cases}, &\forall t \in \{0,...,T-1\}, \\
    & y_{t+1} = \gamma_1 y_t + (1- \gamma_1)y_{b,t} + z_{t+1}, &\forall t \in \{0,...,T-1\}, \\
    & y_{b,t+1} = \gamma_2 y_t + \gamma_3 y_{b,t} + (1- \gamma_2 + \gamma_3)y_b, &\forall t \in \{0,...,T-1\}.
    \end{aligned}
\end{equation}
\end{proposition}
This reformulation is done by using the substitution in the premise of the proposition and adjusting the threshold point in the piece-wise linear dynamics. The full proof can be found in the appendix.
\begin{remark}
Note that since $b>0$, we can recover the solution values of $k,x_t$ in the original optimization problem by dividing $c,z_t$ by the estimated value of $b$ respectively.
\end{remark}
Next we need to consider reformulation the bi-linear terms involving $\alpha z_t$ present in the dynamics.
\begin{proposition}
\label{prop:alpha_z_reform}
For $\alpha_i \in \mathcal{A}$ (for $i = 1,...,m$) let $w_{it} $ and $w_t$ be continuous variables and let  $\iota_{i} \in \mathbb{B}$ be indicator variables. Then Constraints \eqref{eq:hep_dynam_const} can be reformulated as:
\begin{equation} 
    \begin{aligned}
    &z_{t+1} = bu_t + \begin{cases}w_t & \text{for } z_t \leq c + w_t \\ z_t - c & \text{for } z_t > c + w_t \end{cases}, &\forall t \in \{0,...,T-1\}, \\
    & w_t = \sum_{i = 1}^m w_{it}, &\forall t \in \{0,...,T-1\}, \\
    & w_{it} \leq \alpha_i z_t + M(1-\iota_{i}), &\forall t \in \{0,...,T-1\},i \in \{1,..,m\}, \\
    & w_{it} \geq \alpha_i z_t - M(1-\iota_{i}), &\forall t \in \{0,...,T-1\}, i \in \{1,..,m\}, \\
    & \sum_{i=1}^m \iota_i = 1. 
    \end{aligned}
\end{equation}
\end{proposition}
This follows from the standard reformulation for disjunctive constraints \citep{conforti2014integer,wolsey1999integer}. The ``OR'' relationship between different $\alpha_i$ is modeled by indicator variable $\iota_i$. The proof for equivalence between the original problem and reformulated problem can be found in the appendix.

\begin{remark}
The final constraint on $\sum_{i=1}^m \iota_i = 1$ can be modeled with SOS 1 constraint to achieve better solver performance.
\end{remark}
Next we need to reformulate the piece-wise linear dynamics of heparin in such a way that they can be used in a mixed integer linear solver.
\begin{proposition}
\label{prop:MILP_solver}
Let $\nu_t \in \mathbb{B}$ be a set of indicator variables. Then using these variables we can reformulate $g(z_t)$ as:
\begin{equation}
    \begin{aligned}
    & z_{t+1} \leq w_t + bu_t - M\nu_t, &\forall t \in \{0,...,T-1\}, \\
    & z_{t+1} \geq w_t + bu_t + M\nu_t, &\forall t \in \{0,...,T-1\},\\
    & z_{t+1} \leq z_t - c - M(1-\nu_t), &\forall t \in \{0,...,T-1\},\\
    & z_{t+1} \geq z_t - c + M(1-\nu_t), &\forall t \in \{0,...,T-1\}, \\
    & z_t \geq c + w_t - M(1-\nu_t), &\forall t \in \{0,...,T-1\}, \\
    & z_t \leq c + w_t + M\nu_t, &\forall t \in \{0,...,T-1\}. 
    \end{aligned}
\end{equation}
\end{proposition}

This reformulation is done by using big-M constraints to model if-then relationship in the bi-linear terms \citep{wolsey1999integer}. The proof for equivalence between the original problem and reformulated problem can be found in the appendix.

\subsection{Benders Decomposition Training Algorithm}
Although the formulation provided in the previous section can be directly implemented using commercial solvers, it is not practical for real time implementation. We found during implementation that while solvers were able to find a good feasible solution in a reasonable amount of time, they were not able to confirm the solution was optimal within less then 4 hours. This is because, the formulation heavily relies on big-M constraints that can be problematic computationally \citep{bigM}. Due to the structure of our problem, it is difficult to remove all big-M constraints from the final formulation which means that alternative solution methods must be used.

Therefore instead of solving the large-scale optimization problem as a whole we propose using a decomposition approach that can be seen as a special case of a Generalized Benders Decomposition \citep{geoffrion1972generalized} that exploits the structure of optimization problems with complicating variables. When complicating variables are fixed, the remaining optimization problem will be considerably more tractable. To see how this applies to our problem, note that in optimization problem (\ref{eq:constraind_mle}) all integer variables come about from the linearization of the dynamics of $x_t$, which involves variables $\alpha$, $k$ and $b$. If these variables are fixed, we can remove the integer constraints which renders the remaining problem easy to solve. Specifically, consider the following parametric optimization problem for fixed values $\alpha = \bar{\alpha}, k = \bar{k}, b = \bar{b}$.

\begin{equation}
\label{eq:sub_prob}
    \begin{aligned}
\mathcal{L}(\bar{\alpha}, \bar{k},\bar{b}) = \max_{\{z_t,y_t,y_{bt}\}_{t=0}^T,y_b} & \sum_{t \in n(T)} \log p(\tilde{y}_t| y_t) \\
\text{subject to:} &\\
&  z_{t+1} = \bar{b}u_{t+1} + \begin{cases}\bar{\alpha} z_t \quad \text{for } z_t \leq \frac{\bar{b}\bar{k}}{1-\bar{\alpha}} \\ z_t- \bar{b}\bar{k} \quad \text{for } z_t > \frac{\bar{b}\bar{k}}{1-\bar{\alpha}} \end{cases}, & \forall t \in \{0,...,T-1\},\\
& y_{t+1} = \gamma_1 (y_t - y_{bt}) + y_{bt} + z_{t+1}, & \forall t \in \{0,...,T-1\}, \\
&y_{b,t+1} = \gamma_2  y_{bt} + \gamma_3y_{t} + \gamma_4 y_b, & \forall t \in \{0,...,T-1\},\\ 
& \frac{z_t}{\bar{b}} \in \mathcal{X}, y_t,y_{bt},y_b \in \mathcal{Y}, & \forall t \in \{0,...,T\}.
\end{aligned}
\end{equation}

Given this formulation we consider the following proposition:

\begin{proposition}
\label{prop:lp_sub_prob}
For $ \bar{\alpha} \in \mathcal{A},\, \bar{k} \in \mathcal{K},\bar{b} \in \mathcal{B}$ and known control sequence $\{u_t\}_{t=0}^T$, the value function $\mathcal{L}(\bar{\alpha}, \bar{k},\bar{b})$ is the value function of a convex optimization problem with respect to the right hand side of its linear constraints.
\end{proposition}
The key step in deriving this result is to notice that with fixed $\bar{\alpha}, \bar{k},\bar{b}$ and past doses $\{u_t\}_{t=0}^T$, the entire trajectory of heparin $x_t$ can be determined. This means that in (\ref{eq:reform_bilinear}), we are given $\overrightarrow{z}_t = z_t(\bar{\alpha},\bar{k}, \bar{b}, \{u_\tau\}_{\tau=0}^t)$ and the first bilinear term is eliminated. In the remaining constraints, ${z}_t$ appears as an affine term. Therefore, we can reformulate (\ref{eq:sub_prob}) as:

\begin{equation}
\label{eq:lp_sub_z}
    \begin{aligned}
   \mathcal{L}(\overrightarrow{z}(\bar{a},\bar{k},\bar{b})) = \max_{\{z_t,y_t,y_{bt}\}_{t=0}^T,y_b} & \sum_{t \in n(T)} \log p(\tilde{y}_t | y_t) \\
\text{subject to:} &\\
& y_{t+1} = \gamma_1 (y_t - y_{bt}) + y_{bt} + \bar{z}_{t+1}, & \forall t \in \{0,...,T-1\},\\
&y_{b,t+1} = \gamma_2  y_{bt} + \gamma_3y_{t} + \gamma_4 y_b, & \forall t \in \{0,...,T-1\},\\ 
& y_t,y_{bt},y_b \in \mathcal{Y}, & \forall t \in \{0,...,T\}.
    \end{aligned}
\end{equation}

This is a convex optimization problem which can be solved efficiently compared to the original MILP and it yields a lower bound for problem (\ref{eq:constraind_mle}) since we fix a particular combination of $\bar{\alpha}, \bar{k},\bar{b}$ which can be sub-optimal. However, we can now view the original optimization problem as maximizing this parametric function $\mathcal{L}(\overrightarrow{z}(\bar{a},\bar{k},\bar{b}))$ over all feasible combinations of $\bar{a},\bar{k},\bar{b}$. Since the resulting parametric problem is convex and much easier to solve, it has a structure that can be exploited through a Generalized Benders Decomposition approach. To develop the decomposition algorithm, we need to further rewrite problem (\ref{eq:constraind_mle}) into a master problem with respect to the dual problem of $\mathcal{L}(\overrightarrow{z}(\bar{\alpha},\bar{k},\bar{b}))$ and dual feasibility set of $\bar{\alpha}, \bar{k},\bar{b}$. Consider the following proposition:

\begin{proposition}
\label{prop:master}
If Assumptions (1)- (2) hold. The optimization problem in \eqref{eq:constraind_mle} can be reformulated as:
\begin{equation}
\label{eq:master}
    \begin{aligned}
    \max\limits_{\alpha,k,b,\ell}& \  \ell \\
    & \text{subject to:} \\
    & \ell \leq \overrightarrow{z}(\alpha,k,b)^\top \overrightarrow{\lambda} + g(\{\lambda_t\}_{t=0}^{n-1}), & \forall \overrightarrow{\lambda} \in \mathbb{R}^{T-1}, \\
     & 0 = \overrightarrow{z}(\bar{\alpha},\bar{k},\bar{b})^\top \overrightarrow{\mu} + r_y(\{\mu_t\}_{t=0}^{n-1}), & \forall \overrightarrow{\mu} \in [-1,1]^{T-1},  \\
     &\alpha \in \mathcal{A}, k \in \mathcal{K}, b \in \mathcal{B}.
    \end{aligned}
\end{equation}
Where $\overrightarrow{\lambda} = [0,\lambda_0,\lambda_1,....,\lambda_{T-1}]^\top$, $\overrightarrow{\mu} = [0,\mu_0,\mu_1,....,\mu_{T-1}]^\top$ are the vectorization of the sequence of Lagrange multipliers. $g(\{\lambda_t\}_{t=0}^{n-1}) = \sum_{t \in n(T)} \log p(\tilde{y}_t | y^*_t) + \sum_{t=0}^{T-1} \lambda_t\big(\gamma_1(y^*_t -y^*_{bt}) + y^*_{bt} - y^*_{t+1}\big)$ and $y^*_t,y_{bt}^*$ are derived as:
\begin{equation}
\label{eq:dual_sub}
    \begin{aligned}
    \{y^*_t,y_{bt}^*\}_{t=0}^T \in \argmax_{\{y_t,y_{bt}\}_{t=0}^T,y_b }& \sum_{t \in \mathcal{T}_y} \log p(\tilde{y}_t | y_t) + \sum_{t=0}^{n-1} \lambda_t\big(\gamma_1(y_t -y_{bt}) + y_{bt} - y_{t+1}\big) \\
    & \text{subject to:}\\
    &y_{b,t+1} = \gamma_2 y_{bt} + \gamma_3 y_{t} + \gamma_4 y_b, & \forall t \in \{0,...,T-1\}, \\
    &y_t,y_{bt},y_b \in \mathcal{Y}, & \forall t \in \{0,...,T\}.
    \end{aligned}
\end{equation}
while $r_y(\{\mu_t\}_{t=0}^{n-1}) = \sum_{t=0}^{T-1} \lambda_t\big(\gamma_1(\bar{y}^*_t -\bar{y}^*_{bt}) + \bar{y}^*_{bt} - \bar{y}^*_{t+1}\big)$ and $\bar{y}^*_t$, $\bar{y}^*_{bt}$ are derived as:
\begin{equation}
\label{eq:dual_feasi_sub}
    \begin{aligned}
    \{\bar{y}^*_t,\bar{y}_{bt}^*\}_{t=0}^T \in \argmax_{\{y_t,y_{bt}\}_{t=0}^T,y_b }& \sum_{t=0}^{T-1} \mu_t\big(\gamma_1(y_t -y_{bt}) + y_{bt} - y_{t+1}\big) \\
    & \text{subject to:}\\
    &y_{b,t+1} = \gamma_2 y_{bt} + \gamma_3 y_{t} + \gamma_4 y_b, & \forall t \in \{0,...,T-1\}, \\
    &y_t,y_{bt},y_b \in \mathcal{Y}, & \forall t \in \{0,...,T\}.
    \end{aligned}
\end{equation}
\end{proposition}
Notice that (\ref{eq:dual_feasi_sub}) is a linear program and (\ref{eq:dual_sub}) is a convex MILP by our assumption that the noise p.d.f can be described by mixed integer linear constraints but it is still easy to solve in terms of the number of integer constraints it contains. The first constraint in problem (\ref{eq:master}) defines the optimality cuts. It is derived from the Lagrangian dual function for (\ref{eq:lp_sub_z}) with respect to its first constraint: 
\begin{equation}
    \mathcal{L}_{\textit{dual}}(\overrightarrow{z}(\bar{\alpha},\bar{k},\bar{b}))= \min_{\lambda} \overrightarrow{z}(\bar{\alpha},\bar{k},\bar{b})^\top \overrightarrow{\lambda} + g(\{\lambda_t\}_{t=0}^{n-1}).
\end{equation}
This optimization problem is also known as the subproblem in the Benders Decomposition procedure. We want to find $\smash{\displaystyle\max_{\bar{\alpha},\bar{k},\bar{b}}\ \min_{\lambda} \overrightarrow{z}(\bar{\alpha},\bar{k},\bar{b})^\top \overrightarrow{\lambda} + g(\{\lambda_t\}_{t=0}^{n-1})}$. Using the principle of maximin, we get the objective and first constraint in (\ref{eq:master}). The second constraint in (\ref{eq:master}) defines the feasibility cuts that ensure that the feasible set of $\alpha, k, b$ is the same in  (\ref{eq:master}) as  (\ref{eq:constraind_mle}). It is derived from the first constraint in (\ref{eq:lp_sub_z}) since  we are dualizing with respect to the first constraint and we can prove that the second constraint in (\ref{eq:lp_sub_z}) always has a feasible solution. The full proof of these propositions can be found in the appendix.

Problem (\ref{eq:master}) is known as the master problem and it cannot be solved directly because it contains an infinite number of constraints.  A natural way to approach this is relaxation: we start by solving a relaxed version of the master problem that contains only a subset of original constraints. If in the test of feasibility some ignored constraints are violated we add them to the relaxed problem and solve the updated master problem again. The procedure is repeated until a desirable solution of acceptable accuracy has been obtained. With this idea, we propose our $\epsilon$-optimal Benders Decomposition Algorithm as shown in Algorithm \ref{alg:benders_decomp}.

\begin{algorithm}
\caption{$\epsilon$-optimal Benders Decomposition Algorithm}
\label{alg:benders_decomp}
\begin{algorithmic}[1]
\For{$\alpha$ in $\mathcal{A}$}
\State Let a point $(\alpha, \bar{k},\bar{b})$ such that \eqref{eq:feasibility} is feasible given  $\bar{z}_t = z_t(\alpha, \bar{k},\bar{b})$ be known. Solve the Lagrangian dual problem $\mathcal{L}_{dual}(\overrightarrow{z}(\alpha, \bar{k},\bar{b}))$ and obtain the optimal multiplier $\bar{\lambda}$. Let $p=1$, $q = 0$, $\overrightarrow{\lambda}^1=\bar{\lambda}$, $LBD = \mathcal{L}_{dual}(\overrightarrow{z}(\alpha, \bar{k},\bar{b}))$. Select the convergence tolerance parameter $\epsilon > 0$
\While{True}
\State Solve the relaxed master problem and obtain optimal solution $(\ell^*,k^*,b^*)$
\begin{equation}
\label{eq:master_algo}
    \begin{aligned}
    \max_{\ell,\overline{k},\overline{b}} & \  \ell \\
    & \text{subject to:} \\
    & \ell \leq \overrightarrow{z}(\alpha, \bar{k},\bar{b})^\top \overrightarrow{\lambda}^i + g(\{\lambda_t^i\}_{t=0}^{n-1}), & \forall i = 1,...,p,\\
     & 0 = \overrightarrow{z}(\alpha, \bar{k},\bar{b})^\top \overrightarrow{\mu}^i + r_y(\{\mu^i_t\}_{t=0}^{n-1}), & \forall i = 1,...,q.
    \end{aligned}
\end{equation}
\If {$LBD \geq \ell^*-\epsilon$} terminate.
\EndIf
\State Solve subproblem $\mathcal{L}_{dual}(\alpha,k^*,b^*)$ and obtain the optimal multiplier $\overrightarrow{\lambda}$
\If {The quantity $\mathcal{L}_{dual}(\alpha,k^*,b^*)$ is unbounded} Determine  $\overrightarrow{\mu}$ such that $\overrightarrow{z}(\bar{\theta}_x,\bar{b})^\top \overrightarrow{\mu} + r_y(\{\mu_t\}_{t=0}^{n-1}) > 0$. Let $q=q+1$, $\overrightarrow{\mu}^q = \overrightarrow{\mu}$
\Else
\If {$\mathcal{L}_{dual}(\alpha,k^*,b^*)< \ell^* - \epsilon$} obtain the optimal multiplier $\overrightarrow{\lambda}$, let $p=p+1$, $\overrightarrow{\lambda}^p=\overrightarrow{\lambda}$. If $\mathcal{L}_{dual}(\alpha,k^*,b^*) > LBD$, put $LBD = \mathcal{L}(\alpha,k^*,b^*)$
\Else \ Set $\mathcal{L}_{\alpha} = \mathcal{L}_{dual}(\alpha,k^*,b^*)$, $k_{\alpha} = k^*, b_{\alpha} = b^*$ and terminate.
\EndIf
\EndIf
\EndWhile
\EndFor
\State The optimal $\alpha^* = \max_{\alpha \in \mathcal{A}} \mathcal{L}_{\alpha}$ and optimal $k^* = k_{\alpha^*}$, $b^* = b_{\alpha^*}$.
\end{algorithmic}
\end{algorithm}
In Algorithm 1, we iterate between solving the relaxed master problem and the subproblem. The master problem (\ref{eq:master}) yields an upper bound for the optimal value and provides a temporary solution $\alpha^*, k^*, b^*$. This solution might not be feasible because the relaxed master problem does not contain all the constraints. In that case we use the subproblem $\mathcal{L}(\alpha^*,k^*,b^*)$ to test for feasibility If the subproblem is unbounded, it means that $\alpha^*,k^*,b^*$ is not feasible to the original problem. Most dual-type and primal algorithms will yield a $\overrightarrow{\mu}$ with which the dual feasibility constraint is violated. Then we can use $\overrightarrow{\mu}$ in our feasibility cuts in the relaxed master problem. If $\alpha^*,k^*,b^*$ is feasible and the subproblem has a finite optimal solution, then we obtain a lower bound for the original problem. In this case, there are two possibilities: (i) the distance between the lower bound and the upper bound is smaller than $\epsilon$ which means we have arrived at the desired accuracy and the procedure terminates (ii) we add the optimal Lagrangian multiplier $\overrightarrow{\lambda}$  to the optimality cuts in the relaxed master problem. By adding these cuts to the relaxed master problem, we get decreasing upper bounds and increasing lower bounds. In the following proposition, we show the convergence of our $\epsilon$-optimal Benders Decomposition Algorithm.
\begin{proposition}
\label{prop:benders_convergence}
The $\epsilon$-optimal Benders Decomposition Algorithm procedure terminates in a finite number of steps for any given $\epsilon$.
\end{proposition}
This proposition can be obtained by applying canonical result developed in \citep{geoffrion1972generalized} given that the structure our problem matches the assumption listed in the paper. More proof details can be found in the appendix.

The $\epsilon$-optimal Benders Decomposition Algorithm can obtain an approximated solution in a considerably shorter time than solving the original big MILP because it is solving a sequence of much easier optimization problems. The relaxed master problem (\ref{eq:master_algo}) contains fewer constraints defined by different $\overrightarrow{\lambda},\overrightarrow{\mu}$ and the underlying maximization problem $g(\{\lambda_t^i\}_{t=0}^{n-1})$ is not hard to solve. The optimal multiplier to subproblem $\mathcal{L}_{dual}(\alpha,k^*,b^*)$ can be obtained by solving its dual which is a more straight forward optimization problem. In our computational experiments, a near-optimal solution to the original MILP can be obtained within a minute through this decomposition algorithm which makes the patient parameter estimation more applicable in clinical applications.

\subsection{Bayesian Estimation}
\label{sec:bayes_estimation}
The MLE formulation (\ref{eq:mle_prob}) presents one way of obtaining estimators for the patient parameters. However, the clinician may have some prior knowledge about the possible values of these patient parameters, i.e., a prior probability distribution over different combinations of $(\theta_x,\theta_y)$. In this case, a Bayesian framework is a natural setting for making predictions of the patient's aPTT trajectory. Moreover, this approach can be used to quantify the level of uncertainty in the estimation of system state and system dynamic parameters.
Suppose the patient has stayed in the ICU for $T$ time periods with $n(T)$ lab measurements $ \{\tilde{y}_t\}_{t\in n(T)}, \{u_t\}_{t=0}^T$ taken. The clinician wants to predict the patient's condition $\{x_i,y_i\}^{T+n}_{i=t}$ for some $n>0$ steps into the future. In a Bayesian view, this means the clinician wants to calculate the posterior distribution of $\{x_i,y_i\}^{T+n}_{i=t}$. As $(\theta_x,\theta_y)$ completely characterize the patient dynamics, this task can be done by calculating the posterior distribution of $(\theta_x,\theta_y)$. A direct application of Bayes's Theorem \citep{bickel2015mathematical} shows that the posterior of of $(\theta_x,\theta_y)$ can be expressed as
\begin{equation}
    p\big(\theta_x,\theta_y  | \{\tilde{y}_t\}_{t\in n(T)}, \{u_t\}_{t=0}^T \big) = Z^{-1} \times p\big(\{\tilde{y}_t\}_{t\in n(T)} | \theta_x,\theta_y, \{u_t\}_{t=0}^T \big) \times p(\theta_x,\theta_y).
\end{equation}
Here $Z$ is a normalization constant that ensures the right hand side is a probability distribution and $p(\theta_x,\theta_y)$ represents the clinician's prior knowledge. To ensure the resulting optimization problem is solvable by commercial solvers, we need an assumption on $p(\theta_x,\theta_y)$.
\begin{assumption}
\label{post_assumption}
The function $\log p(\theta_x,\theta_y)$ can either be expressed using a finite number of mixed integer linear constraints and is concave, and $p(\theta_x,\theta_y) > 0$ for all $(\theta_x,\theta_y) \in \Theta_x \times \Theta_y$.
\end{assumption}
This is a mild assumption because it holds for the Laplace distribution, the shifted exponential distribution, and piece-wise linear distributions. This ensures that all possible parameter values will have non-zero probability density. Moreover this is satisfied by a uniform prior distribution, which is equivalent to solving the MLE problem. Next, we introduce a profile likelihood \citep{murphy2000profile} approach to computing the posterior distribution of $(\theta_x,\theta_y)$.
First, consider the following problem:
\begin{equation}
\label{eq:posterior_optimization}
    \begin{aligned}
    \psi_{T}(\bar{\theta}_x,\bar{\theta}_y) &= \log p\big(\theta_x = \bar{\theta}_x ,\theta_y = \bar{\theta}_y | \{\tilde{y}_t\}_{t\in n(T)}, \{u_t\}_{t=0}^T \big) + \log Z,\\
    =\max& \sum_{t\in n(T)} \log p(\tilde{y}_t|y_t) + \log p(\theta_x, \theta_y)\\
\text{subject to: } &\\
& \eqref{eq:hep_dynam_const}-\eqref{eq:ptt_baseline_const}, &  \forall t \in \{0,...,T-1\}, \\
& x_0 = \bar{x}_0, y_0 = \bar{y}_0, y_{b,1} = \bar{y}_{b,1}, y_b = \bar{y}_b, \\
& x_t \in \mathcal{X}, y_t,y_{bt},y_b, \in \mathcal{Y} & \forall t \in \{0,...,T\}.\\ 
    \end{aligned}
\end{equation}
Notice that the above problem only differs from the MLE formulation (\ref{eq:mle_prob}) in two ways: first, the initial states and the patient parameters for the dynamics are fixed; second, there is an additional term in the objective $\log p(\theta_x,\theta_y)$. Therefore, with reformulation introduced in Section \ref{sec:milp_form_mle} and Assumption \ref{post_assumption}, the above problem (\ref{eq:posterior_optimization}) can be expressed as a MIP. Solving (\ref{eq:posterior_optimization}) does not directly provide the posterior distribution of $(\theta_x,\theta_y)$ because $Z$ is not known a \emph{priori}. But since $Z$ only scales the posterior estimate, we instead use a simpler scaling similar to that proposed by \cite{mintz2017behavioral}. Let $(\hat{\theta}_x,\hat{\theta}_y) \in \argmax_{(\theta_x,\theta_y)}\psi_{T}(\theta_x,\theta_y)$ be the maximum a \emph{posteriori} estimates of the patient parameters and note that these estimates can be obtained by solving (\ref{eq:posterior_optimization}) with fixed parameter constraints removed which is a MIP. We propose using:
\begin{equation}
\label{bayesian_estimator}
    \hat{p}\big(\theta_x,\theta_y  | \{\tilde{y}_t\}_{t\in n(T)}, \{u_t\}_{t=0}^T \big) = \frac{\exp(\psi_{T}(\theta_x,\theta_y))}{\exp(\psi_{T}(\hat{\theta}_x,\hat{\theta}_y))}.
\end{equation}
as an estimate of the posterior distribution of $(\theta_x,\theta_y)$. Two key properties of our new estimate are that $\hat{p}\big(\theta_x,\theta_y  | \{\tilde{y}_t\}_{t\in n(T)}, \{u_t\}_{t=0}^T \big) \in [0,1]$ by construction, and that $\hat{p}\big(\hat{\theta}_x,\hat{\theta}_y  | \{\tilde{y}_t\}_{t\in n(T)}, \{u_t\}_{t=0}^T \big) = 1$.

\subsection{Consistency of Estimates}
In this section we will discuss the statistical properties of our estimates. Mainly, we show that the estimate (\ref{bayesian_estimator}) is consistent in a Bayesian sense \citep{bickel2015mathematical}.

\begin{definition}
The posterior estimate \eqref{eq:posterior_optimization} is consistent if for all $(\theta^*_x,\theta^*_y) \in \Theta_x \times \Theta_y $and $\delta, \epsilon> 0$ we have $p_{(\theta^*_x,\theta^*_y)}(\hat{p}(\mathcal{E}(\delta))|\{\tilde{y}_t\}_{t\in n(T)}, \{u_t\}_{t=0}^T)\geq \epsilon) \rightarrow 0$ as $T \rightarrow 0$ where $p_{(\theta^*_x,\theta^*_y)}$ is the probability law under $(\theta^*_x,\theta^*_y)$, $\mathcal{E}(\delta)=\left\{\theta_x,\theta_y) \notin \mathcal{B}\left(\theta^*_x,\theta^*_y, \delta\right)\right\}$ where $\mathcal{B}\left(\theta^*_x,\theta^*_y, \delta\right)$ is an open $\delta$ ball around $(\theta^*_x,\theta^*_y)$.
\end{definition}

This definitions means that for any point other than the true initial parameters $(\theta^*_x,\theta^*_y)$, the log likelihood becomes infinitely small as we have more observations. We need an additional assumption known as \emph{sufficient excitation} to prove the statistical consistency of \eqref{eq:mle_prob}.

\begin{assumption}
\label{sufficient_excitation}
Let $(\theta^*_x,\theta^*_y)$ be the patient's true parameters. The heparin doses $u_t$ are such that

\begin{equation}
    \max _{\mathcal{E}(\delta)} \lim _{T \rightarrow \infty} \sum_{t\in n(T)} \log \frac{p_{\epsilon}\left(\tilde{y}_{t}- \bar{y}_{t}\right)}{p_{\epsilon}\left(\tilde{y}_{t}- y^*_{t}\right)}=-\infty,
\end{equation}

for any $\delta > 0$, almost surely, where $y^*_t$ are states under true parameters $(\theta^*_x,\theta^*_y)$, and $\bar{y}_t$ are the states under any other parameters $(\theta_x,\theta_y)$.
\end{assumption}
This type of assumption is common in the adaptive control literature \citep{craig1987adaptive, aastrom2013adaptive} known as a \emph{sufficient excitation} or a \emph{sufficient richness} condition. Heparin dosing conditions in practice satisfy this assumption due to the noise in aPTT measurements.
\begin{proposition}
\label{corol:bayesian_consistency}
If Assumptions \ref{as:alpha_set_assump}-\ref{sufficient_excitation} hold then the posterior estimate $ \hat{p}\big(\theta_x,\theta_y  | \{\tilde{y}_t\}_{t\in n(T)}, \{u_t\}_{t=0}^T \big)$ given by \eqref{bayesian_estimator}  is consistent.
\end{proposition}
To prove this result, we expand the expression for $\log  \hat{p}\big(\theta_x,\theta_y  | \{\tilde{y}_t\}_{t\in n(T)}, \{u_t\}_{t=0}^T \big)$ and write it in terms of log probability of observation error $p_{\epsilon}\left(\tilde{y}_{t}- {y}_{t}\right)$ and log transition probability $p(x_{t-1}, u_t, \theta_x|x_t)$. Then we use the fact that that transition probability $p(x_{t-1}, u_t, \theta_x|x_t)$ is degenerate to show $\log \hat{p}(\theta_x, \theta_y)$ diverges to $-\infty$ for any parameter other than the ground truth. Then the uniform result on the $\delta$ ball can be derived using the volume bound. The full proof details can be found in the appendix.

\begin{corollary}
\label{corol:estimator_consistency}
If Assumptions \ref{as:alpha_set_assump}-\ref{sufficient_excitation} hold, then the maximum a \emph{posteriori} estimates $(\hat{\theta}_x, \hat{\theta}_y) \underset{}{\stackrel{p}{\longrightarrow}} (\theta^*_x,\theta^*_y)$ as $T \rightarrow \infty$.
\end{corollary}

This result can be proved by showing that with probability 1, $\hat{\theta}_x, \hat{\theta}_y$ is included in a $\delta$-ball around $\theta^*_x, \theta^*_y$ as $T \rightarrow \infty$ which is implied by Proposition \ref{corol:bayesian_consistency}. The full proof is in the appendix. 

The above two corollaries imply that the MAP and MLE estimators are consistent. Also note that as stated in Section \ref{sec:bayes_estimation}, these estimators can be calculated using Algorithm \ref{alg:benders_decomp}.

\section{Dose Optimization}
\label{sec:dose_opt}
In the previous section, we developed a methodology that provides consistent estimates of model parameters that characterize the system state and system dynamics for each patient. Knowing a patient's individual parameters gives us predictive insight into how the patient's heparin and aPTT will evolve in response to any dosing schedule. In this section, we will leverage this patient-specific model to decide future doses that can keep the patient's aPTT level in a safe therapeutic range. Specifically, we consider the following adaptive framework for  dose design: first we estimate the patient's parameters using their aPTT lab measures, then using the estimates and the uncertainty of the estimation we solve an optimization problem to design the optimal dose. We will repeat these steps every 4-6 hours as new measurements of aPTT are collected, to both improve the estimation accuracy and adaptively change the dose amounts to address the patient's current condition.
%
%
One of the key challenges of this problem is that improper estimation of patient parameters can lead to misdosing that can cause adverse health effects. 
To address this safety concern, we will provide theoretical guarantees that the dosing sequences output by our method are appropriate for the patient. In this section, we develop a scenario generation based dosing algorithm that outputs asymptotically optimal dosing sequences with respect to the patient's true parameters.

\subsection{Problem Formulation and Preliminaries}
\label{dose_formulation}
To formally formulate our problem, let $\ell: \mathcal{Y}^n \times \mathcal{U}^n \rightarrow \mathbb{R}$ be a bounded loss function of the patient's aPTT level and dose over the next $n$ time points that reflects the need for aPTT to follow a desired therapeutic trajectory.  If the current time period is time $T$, our goal will be to find a dosing sequence $\{u_t\}_{t=T}^{T+n} \in \mathcal{U}^n$ for each patient such that the loss is minimized. For our analysis we make the following fairly general assumption on the structure of $\ell$.


\begin{assumption}
\label{loss_function}
The loss function $\ell$ can be described by mixed integer linear constraints.
\end{assumption}
This assumption is not restrictive since it applies to large class of piece-wise linear functions. We provide concrete examples of potential loss functions in our experiments in Section \ref{experiments}.

Since the clinician only has noisy and incomplete observations of the patients aPTT level $\{\tilde{y}_t\}_{t\in n(T)}$, one planning approach is to minimize the expected posterior loss associated with applying certain dosing policy $\left\{u_{t}\right\}_{t=T+1}^{T+n}$ in next $n$ periods:
\begin{equation}
    \min_{\{u_t\}_{t=T+1}^{T+n}} \left\{\mathbb{E}\left[\ell \left(\left\{y_{t}, u_{t}\right\}_{t=T+1}^{T+n}\right)\mid \left\{\tilde{y}_{t}\right\}_{t\in n(T)},\left\{u_{t}\right\}_{t=0}^{T} \right] \bigg| \left\{u_{t}\right\}_{t=T+1}^{T+n} \in \mathcal{U}^n\right\}.
\end{equation}

Recall that the patient's state trajectory is completely characterized by the parameters $(\theta_x, \theta_y)$, and so by the sufficiency and the smoothing theorem \citep{bickel2015mathematical}, there exists $\varphi: \Theta_x \times \Theta_y \times \mathcal{U}^n \mapsto \mathbb{R}$ such that the planning problem can be formulated as:
\begin{equation}
  \min_{\{u_t\}_{t=T+1}^{T+n}} \left\{\mathbb{E}\left[\varphi \left(\theta_x, \theta_y, \{u_{t}\}_{t=0}^{T+n}\right) \mid \left\{\tilde{y}_{t}\right\}_{t\in n(T)},\left\{u_{t}\right\}_{t=0}^{T}\right] \bigg|\left\{u_{t}\right\}_{t=T+1}^{T+n} \in \mathcal{U}^n\right\}.  
\end{equation}

In practice we do not know the true parameters $(\theta^*_x, \theta^*_y)$ so our goal is to design a good approximation to the function $\varphi \left(\theta^*_x, \theta^*_y, \{u_{t}\}_{t=0}^{T+n}\right)$ using our estimation framework. One choice would be simply using the MLE parameters as plug in estimators $\hat{\theta}_{x,MLE}, \hat{\theta}_{y,MLE}$ obtained by solving \eqref{eq:mle_prob} and hoping that $\varphi \left(\hat{\theta}_{x,MLE}, \hat{\theta}_{y,MLE}, \{u_{t}\}_{t=0}^{T+n}\right)$ is convergent to $\varphi \left(\theta^*_x, \theta^*_y, \{u_{t}\}_{t=0}^{T+n}\right)$. However, this approach does not fully address the challenge of safety since the estimates $\hat{\theta}_{x,MLE}, \hat{\theta}_{y,MLE}$ can have a large amount of estimation variance. Also, in general point-wise convergence of a sequence of stochastic optimization problems is not sufficient to ensure convergence of the minimizers of the sequence of optimization problems to the minimizer of the limiting optimization problem \citep{rockafellar2009variational}. In fact, experimental results confirm that this plug-in does not perform well in our setting. Therefore, we propose using a scenario generation approach to addressing the approximation of function $\varphi \left(\theta^*_x, \theta^*_y, \{u_{t}\}_{t=0}^{T+n}\right)$ such that the convergence property of its minimizers is guaranteed. For the purposes of our analysis we introduce the following assumption:
\begin{assumption}
\label{b_assumption}
The set $\mathcal{B}$ for possible values of $b$ is finite, that is $|\mathcal{B}| < \infty$.
\end{assumption}
While our framework discussed in Section \ref{sec:m_param_est} works perfectly when $\mathcal{B}$ is a compact set, this new assumption is necessary for the development and theoretical analysis of our dosing algorithm. Concretely, if we discretize $b$ and fix it in the optimization problem \eqref{eq:constraind_mle}, the objective function will obtain lower semicontinity with respect to the remaining variables which is a key property used in the analysis of optimizers. In practice, the discretization of the original set $\mathcal{B}$ will not significantly impact the accuracy of parameter estimation because we can always approximate the compact set using a $\epsilon$ grid in $\mathcal{B}$ for arbitrarily small $\epsilon$. The discretization allows us to consider an optimization problem that involves a subset of model parameters:
\begin{equation}
    \psi'_T(\bar{\alpha},\bar{b}) = \max_{k,y_b,y_0,y_{b0}} \log\hat{p}(k,y_b,y_0,y_{b0}|\bar{\alpha},\bar{b},\{\tilde{y}_t\}_{t\in n(T)}, \{u_t\}_{t=0}^T).
\end{equation}
Notice that this optimization problem is still a MILP since it is essentially problem (\ref{eq:posterior_optimization}) with $\bar{\alpha}, \bar{b}$ fixed in constraints and $k,y_b,y_0,y_{b0}$ as decision variables. Also denote the optimizers of the above problem as:
\begin{equation}
    \tau(\bar{\alpha},\bar{b}) = \argmax_{k,y_b,y_0,y_{b0}} \hat{p}(k,y_b,y_0,y_{b0}|\bar{\alpha},\bar{b},\{\tilde{y}_t\}_{t\in n(T)}, \{u_t\}_{t=0}^T).
\end{equation}
For fixed observation sequence $\{\tilde{y}_t\}_{t\in n(T)}$ and past doses $ \{u_t\}_{t=0}^T$, $\tau(\bar{\alpha},\bar{b})$ can be interpreted as the profile likelihood estimate  \citep{murphy2000profile} of parameters $k,y_b,y_0,y_{b0}$ with $\alpha, b$ profiled out.
Next we consider an optimization problem parameterized by $\bar{\alpha},\bar{b}$:
\begin{equation}
\label{eq:optimal_dosing}
    \begin{aligned}
    \varphi_{\bar{\alpha},\bar{b}}\left(k,y_b,y_0,y_{b0},\{u_t\}_{t=0}^{T+n}\right) &= \min_{k,y_b,y_0,y_{b0},\{u_t\}_{t=0}^{T+n}} \sum_{t=T+1}^{T+n} \ell(y_t,u_t)\\
\text{subject to:} &\\
&  x_{t+1} = u_{t+1} + \begin{cases}\bar{\alpha} x_t \quad \text{for } x_t \leq \frac{k}{1-\bar{\alpha}} \\ x_t-k \quad \text{for } x_t >\frac{k}{1-\bar{\alpha}} \end{cases}, & \forall t \in \{0,...,T+n\},\\
& y_{t+1} = \gamma_1 (y_t - y_{bt}) + y_{bt} + \bar{b}x_{t+1}, & \forall t \in \{0,...,T+n\}, \\
&y_{b,t+1} = \gamma_2y_{bt} + \gamma_3 y_{t} + \gamma_4 y_b, & \forall t \in \{0,...,T+n\},\\
& x_0=0,x_t \in \mathcal{X}, y_t,y_{bt},y_b \in \mathcal{Y}, & t \in \{1,...,T+n\}.\\ 
    \end{aligned}
\end{equation}

Here $\varphi_{\bar{\alpha},\bar{b}}\left(k,y_b,y_0,y_{b0},\{u_t\}_{t=0}^{T+n}\right)$ is the value function \citep{ralphs2014value} of a parameteric optimization problem with parameters  $k,y_b,y_0,y_{b0},\{u_t\}_{t=0}^{T+n}$ that all belong to affine terms in the constraints. Notice that this is also a feasibility problem since all the parameters are given and dose sequence up to time $T+n$ is fixed. Therefore by Assumption \ref{loss_function} this is a MILP and all the reformulation applied to the MLE problem \eqref{eq:constraind_mle} can be applied here.

\subsection{Adaptive Dosing Algorithm}

Here, we present our adaptive dosing algorithm called Predict then Control with $m$ Scenario Generation (PTC-SGm) in Algorithm \ref{alg:scenario_generation}.

\begin{algorithm}
\caption{Predict then Control with $m$ Scenario Generation (PTC-SGm)}
\label{alg:scenario_generation}
\begin{algorithmic}[1]
\Require $\{\tilde{y}_t\}_{t\in n(T)}$, $\{u_t\}^{T}_{t=0},  \text{ $m$ scenarios } |\bar{\mathcal{A}} \times \bar{\mathcal{B}}| = m$
\State Compute $\omega_{\bar{\alpha},\bar{b}} = \exp\psi'_T(\bar{\alpha},\bar{b})$ for all $\bar{\alpha},\bar{b} \in \bar{\mathcal{A}} \times \bar{\mathcal{B}}$ \newline
\hspace*{4em}$\psi'_T(\bar{\alpha},\bar{b}) = \max_{k,y_b,y_0,y_{b0}} \log\hat{p}(k,y_b,y_0,y_{b0}|\bar{\alpha},\bar{b},\{\tilde{y}_t\}_{t\in n(T)}, \{u_t\}_{t=0}^T)$ \newline
\hspace*{4em} $\tau(\bar{\alpha},\bar{b}) = \argmax_{k,y_b,y_0,y_{b0}} \hat{p}(k,y_b,y_0,y_{b0}|\bar{\alpha},\bar{b},\{\tilde{y}_t\}_{t\in n(T)}, \{u_t\}_{t=0}^T)$
\State Obtain $u_{PTC-SGm}(T) = \{u_t^*\}_{t=T+1}^{T+n} = \argmin_{\{ u_t\}_{t=T}^{T+n}} \frac{\sum\limits_{\bar{\alpha},\bar{b} \in \bar{\mathcal{A}} \times \bar{\mathcal{B}}} \varphi_{\bar{\alpha},\bar{b}}(\tau(\bar{\alpha},\bar{b}),\{u_t\}_{t=0}^{T+n})\omega_{\bar{\alpha},\bar{b}}}{\sum\limits_{\bar{\alpha},\bar{b} \in \bar{\mathcal{A}} \times \bar{\mathcal{B}}} \omega_{\bar{\alpha},\bar{b}}}$
\end{algorithmic}
\end{algorithm}

The PTC-SGm algorithm is inspired by scenario generation method \citep{scenario} for stochastic programming which is a common approach to discretization. Instead of relying on one set of patient parameters obtained by MLE, we assume that there is a distribution over different combinations of $\alpha$ and $b$ and the distribution is determined by
profile likelihood of each combination. Thus, each potential value of $\alpha,b$ can be thought of as a different scenario that we would need to account for. PTC-SGm has two phases that are executed at each time $T$. First, we calculate the distribution of $\alpha, b$, that is we compute the likelihood of each combination $\alpha, b$ in $\mathcal{A}\times\mathcal{B}$ as $\omega_{\alpha,b}$. This step can be thought of as both predicting future states and quantifying the uncertainty around the parameter estimates. Then, we decide doses by solving a problem where each $\alpha,b$ scenario's term in the objective is weighed by $\omega_{\alpha,b}$, while regardless of scenario control sequence $\{u_t\}_{t=0}^{T+n}$ must satisfy appropriate constraints. At each iteration only dose $u_{T}$ will be deployed to the patient, and the algorithm steps will be repeated as new measurements are received. Intuitively, the PTC-SGm algorithm iterates between weighing different combinations of parameters properly and optimizing doses for next period based on weights for different scenarios.

\subsection{Asymptotic Optimality}
An important property of this two-step approach is that the PTC-SGm algorithm outputs asymptotically optimal doses in line with true patient parameters. We break the proof of this result into several steps. The first step, Proposition \ref{algorithm_asymptotic1} builds the relationship between the weighted objective function in the control stage in PTC-SGm and the objective function in optimal dosing problem (\ref{eq:optimal_dosing}) with true patient parameters. 

\begin{corollary}
\label{coro:dose}
If Assumptions \ref{as:alpha_set_assump}-\ref{loss_function} hold, then $\varphi_{\bar{\alpha},\bar{b}}\left(k,y_b,y_0,y_{b0},\{u_t\}_{t=0}^{T+n}\right)$ is lower semicontinuous in $k,y_b,y_0,y_{b0},\{u_t\}_{t=0}^{T+n}$.
\end{corollary}
The proof of this corollary relies on the reformulation introduced in Section \ref{sec:milp_form_mle}. After the reformulation, all the decision variables belong to an affine term in the mixed integer linear constraints. Using the standard results \citep{ralphs2014value} from optimization theory the corollary follows. Proposition \ref{algorithm_asymptotic1} builds the relationship between the weighted objective function in the control stage in PTC-SGm and the objective function in optimal dosing problem (\ref{eq:optimal_dosing}) with true patient parameters. The weighted objective function is a reasonably good approximation to the true optimal dosing objective because it converges to the true objective in probability.

\begin{proposition}
\label{algorithm_asymptotic1}
Suppose that Assumptions \ref{as:alpha_set_assump}-\ref{b_assumption} hold. Then as $T \to \infty$ and for all fixed $\{u_t\}_{t=0}^{T+n}$ we have that $\frac{\sum\limits_{\bar{\alpha},\bar{b} \in\bar{\mathcal{A}} \times \bar{\mathcal{B}}} \varphi_{\bar{\alpha},\bar{b}}(\tau(\bar{\alpha},\bar{b}),\{u_t\}_{t=0}^{T+n})\omega_{\bar{\alpha},\bar{b}}}{\sum\limits_{\bar{\alpha},\bar{b} \in \bar{\mathcal{A}} \times \bar{\mathcal{B}}} \omega_{\bar{\alpha},\bar{b}}} \stackrel{p}{\longrightarrow}  \varphi_{\alpha^*,b^*}\left(k^*,y_b^*,y_0^*,y_{b0}^*,\{u_t\}_{t=0}^{T+n}\right)$.
\end{proposition}
The proof of this proposition has two key parts. The first part is to show that the probability measure $\frac{\omega_{\alpha,b}}{\sum\limits_{\bar{\alpha},\bar{b} \in \bar{\mathcal{A}} \times \bar{\mathcal{B}}} \omega_{\bar{\alpha},\bar{b}}}$ is degenerate at the true parameters as $T$ approaches infinity using Corollary \ref{corol:estimator_consistency}. The second part relies on the property of profile likelihood estimation to show that the optimization problem $\varphi_{\alpha^*,b^*}(\tau(\bar{\alpha},\bar{b}),\{u_t\}_{t=0}^{T+n})$ converges to the one with respect to true parameters in probability. Combing these two limiting results together we complete the proof.

By Proposition \ref{algorithm_asymptotic1} we get a point-wise convergence result. To show the convergence of the optimizers, we need the function $\frac{\sum\limits_{\bar{\alpha},\bar{b} \in\bar{\mathcal{A}} \times \bar{\mathcal{B}}} \varphi_{\bar{\alpha},\bar{b}}(\tau(\bar{\alpha},\bar{b}),\{u_t\}_{t=0}^{T+n})\omega_{\bar{\alpha},\bar{b}}}{\sum\limits_{\bar{\alpha},\bar{b} \in \bar{\mathcal{A}} \times \bar{\mathcal{B}}} \omega_{\bar{\alpha},\bar{b}}}$ to be a lower semicontinuous approximation \citep{vogel2003continuous} to $\varphi_{\alpha^*,b^*}\left(k^*,y_b^*,y_0^*,y_{b0}^*,\{u_t\}_{t=0}^{T+n}\right)$ uniformly in $\mathcal{U}^n$. 
\begin{proposition}
\label{algorithm_asymptotic2}
Suppose that Assumptions \ref{as:alpha_set_assump}-\ref{b_assumption} hold. Then as $T \to \infty$ and for all fixed $\{u_t\}_{t=0}^{T+n}$ we have that $\frac{\sum\limits_{\bar{\alpha},\bar{b} \in \mathcal{A} \times \mathcal{B}} \varphi_{\bar{\alpha},\bar{b}}(\tau(\bar{\alpha},\bar{b}),\{u_t\}_{t=0}^{T+n})\omega_{\bar{\alpha},\bar{b}}}{\sum\limits_{\bar{\alpha},\bar{b} \in \mathcal{A} \times \mathcal{B}} \omega_{\bar{\alpha},\bar{b}}} \underset{\mathcal{U}^n}{\stackrel{l-p r o b}{\longrightarrow}} \varphi_{\alpha^*,b^*}\left(k^*,y_b^*,y_0^*,y_{b0}^*,\{u_t\}_{t=0}^{T+n}\right))$. Here  $\Lambda_n \underset{\mathcal{X}}{\stackrel{l-p r o b}{\longrightarrow}} \Lambda$ means random function $\Lambda_n: \mathcal{X} \rightarrow \mathbb{R}$ is a lower semicontinuous approximation to function $\Lambda: \mathcal{X} \rightarrow \mathbb{R}$.
\end{proposition}

To prove this proposition, we first use the relationship between point-wise lower semicontinuous approximation and convergence in probability. Then we use the fact that point-wise lower semicontinuity everywhere on the space implies a uniform lower semicontinuity for value functions of minimization problems \citep{vogel2003continuous} to conclude the result.

\begin{theorem}
\label{asymptotic_doses}
Note that $\argmin \{\varphi_{\alpha^*,b^*}\left(k^*,y_b^*,y_0^*,y_{b0}^*,\{u_t\}_{t=0}^{T+n}\right)|\{u_t\}_{t=T+1}^{T+n} \in \mathcal{U}^n\}$ is the set of optimal dose sequences under the patient's true parameters $(\theta^*_x,\theta^*_y)$. If Assumptions \ref{as:alpha_set_assump}-\ref{b_assumption} hold, then:
\begin{equation}
    \operatorname{dist}\left(u_{PTC-SGm}(T), \argmin \{\varphi_{\alpha^*,b^*}\left(k^*,y_b^*,y_0^*,y_{b0}^*,\{u_t\}_{t=0}^{T+n}\right)|\{u_t\}_{t=T+1}^{T+n} \in \mathcal{U}^n\}\right) \stackrel{p}{\rightarrow} 0
\end{equation}
as $T \rightarrow \infty$ for any $u_{PTC-SGm}(T)$ returned by PTC-SGm. Here $\operatorname{dist}(x, B)=\inf _{y \in B}\|x-y\|$.
\end{theorem}
With the above two propositions, this theorem can be directly obtained by applying results from stochastic programming theory \citep{vogel2003continuous}.
This result states that the optimal dosing sequence returned by PTC-SGm  is asymptotically included within the set of optimal doses computed with the patient's true parameters, i.e, the dosing sequence obtained by PTC-SGm algorithm is asymptotically optimal. Note that this is not trivial nor is it obvious. The property of lower semi-continuity is key since it ensures the convergence of optima of stochastic optimization problems as a result of the consistency of our posterior estimate. Without these results we would not have the guarantee that our computed policies improve as additional observations are collected. The full proofs of all these results can be found in the appendix.

\section{Computational Results}
\label{experiments}

In this section we discuss two sets of computational studies that evaluate our heparin patient model and the PTC-SGm algorithm. The first set focuses on the predictive accuracy of the dynamic heparin model. We fit the model to the observed aPTT trajectories of patients in a data set that contains real ICU data and predict their future aPTT values under a given sequence of doses from the data set. We then compare the predictive performance of our model against existing machine learning approaches. The second set focuses on examining the performance of PTG-SGm. We simulate a cohort of patients receiving heparin treatment according to variations on PTG-SGm, a naive policy, and existing weight based protocols. The simulation is constructed using the same data set as the predictive experiments.
We evaluate the performance of our algorithm using several metrics that capture treatment efficacy, safety, and tractability. 
%

\subsection{MIMIC III Data}

Our computational experiments are based on a publicly available data set called Medical Information Mart for Intensive Care (MIMIC) III \citep{mimic3} which contains data associated with more than 50,000 distinct hospital admissions for adult patients admitted to ICUs between 2001 and 2021. The original data set is comprehensive and includes various type of data such as patient demographics, laboratory test results, observation notes, and more. For the purpose of our experiments, we extracted a subset of data of patients who have at least 20 records of heparin doses in their chart data and aPTT lab measurements. We used the chart and lab records from 24 hours prior to starting heparin treatment until 24 after ending heparin treatment of a remaining cohort of 25 patients.

\subsection{Predictive Model Comparison}
\label{pre_model_com}
In this computational study, we examine the predictive performance of our heparin patient model in predicting future aPTT values for individual patients. We compare our model against other existing heparin dosing and aPTT prediction methods in the literature that include various state-of-the-art machine learning methods. Since many of the existing methods were developed for classification and not regression we need to convert the prediction task of estimating aPTT into a classification problem. To do this, using the MLE estimates of each patient's baseline aPTT $y_b$ we compute their individual therapeutic range as $[1.5y_b,2.5y_b]$ \citep{basu1972prospective}. Using this range we then label each patient at each time period as either in therapeutic range, below the therapeutic range (i.e. sub-therapeutic), or above the therapeutic range (i.e. super-therapeutic). To reflect a prediction setting that is similar to the deployment of dosing schedules, we use aPTT observations and aggregate heparin data in 4 hour intervals. Our prediction methods then use all historic data up to the time of prediction to predict the aPTT label in the beginning of the next 4 hour epoch. 

\subsubsection{Comparative Models}
We compare our dynamic model against several machine learning methods ranging from regression models to neural networks with more complex architecture. In particular the version of our dynamic patient model we used for the comparison was only the MLE implementation of parameter estimation (i.e. assuming a uniform prior over model parameters). We then converted its regression outputs into a binary prediction using a similar method as proposed in \cite{aswani2019behavioral} by using numerical integration. For the parameter space of the MLE, we specified that $\mathcal{A}=\{0.500, 0.574, 0.630, 0.673, 0.707\}$ which corresponds to the half lives of heparin from 30 minutes to 2 hours in 15 minute increments, $\mathcal{B}=[0.1,10], \mathcal{K}=[0.1,5000], \mathcal{Y}=[0,150]$, and $\mathcal{X}=[0,5000]$.  For comparison, we considered the following classical machine learning methods:

\textbf{Logistic Regression}
We fit two logistic regression models to predict the probability of the patient being in sub-therapeutic range and in super-therapeutic range respectively and decide the patient's status. We use both medical and demographic features as dependent variables. These features include patient's past aPTT observations, time spent in the ICU, past heparin doses, hospital admission type, gender, age, ethnicity, weight and height. The model implementation is made to resemble the model deployed in \cite{multilogit}.

\textbf{Multinomial Logistic Regression}
Multinomial logistic regression is a generalization of logistic regression that allows for direct multi-class classification. Essentially, instead of assuming a Bernoulli distribution over the labels the model now assumes a categorical distribution. With a single model it is designed to address multi-class classification tasks  and works best when the the categories are nominal \citep{greene2003econometric}. We use the same dependent variables in the multinomial logistic regression as for logistic regression.

We also consider deep learning prediction methods and various Artificial Neural Network (ANN) architectures and we use the same features mentioned above as inputs to the neural nets. These methods are calibrated to represent those found in \cite{deeparch} and \cite{Nemati2016OptimalMD}, and have their parameters chosen in  a similar manner to these papers. 

\textbf{Feed Forward ANN}
We construct a feed forward ANN with 1 hidden layer of 10 hidden units. Each  unit in the hidden layer has ReLU activation function and for the output layer we use softmax activation \citep{hagan1997neural}. This method can be seen as a further generalization of the multinomial logistic regression, with the addition of the hidden units allowing for capturing more complex dependencies between features. We used the same features for the feed forward ANN as the the classical machine learning methods.

\textbf{Long short term-memory (LSTM)}
LSTMs are a special kind of recurrent neural network (RNN) capable of capturing long-term dependencies in data streams \citep{hochreiter1997long}. LSTM's maintain an estimate of an internal latent state and allows for prediction of sequences of labels and not just single labels. This makes them a strong candidate for the prediction of aPTT. Our LSTM network included 4 sequential LSTM layers with ReLu activation and 16 hidden units and a final feed forward layer with three outputs and softmax activation. The inclusion of 4 LSTM units is meant to model the 4 hour period between observations of aPTT.

\textbf{Gated Recurrent Units (GRU)}
GRUs are another variation of RNN. They are designed to mitigate the vanishing gradient problem encountered in RNNs \citep{cho2014learning}. As a type of RNN, these models are also capable of predicting sequences of labels and not only single labels. We used a similar architecture for the GRU model as we did for the LSTM model with the key difference being the replacement of the LSTM units with GRU units.
%

The logistic regression and multinomial logistic regression models were trained using the Logistic regression package of scikit learn \citep{scikit-learn} with standard cross validation parameters, while the remaining models were trained using PyTorch \citep{paszke2017automatic}.

For the feed forward ANN, logistic regression, and multinomial  logistic regression we used leave one out cross validation to train the parameters of the models and evaluate predictive performance. The GRU and LSTM models are more complex to train then the classical machine learning models and feed forward ANN. For these models, we still used leave one out cross validation, however the data used in training was transformed into a form containing only one step transitions between data points. This formed the basis of an experience replay memory from which we re-sampled sequences of state transitions of varying lengths. The gradient was then computed along the entirety of each sequence before updating the parameters of the GRU and LSTM models. Due to the dependence of these models on sufficient one step transitions (and thus not much missing aPTT data), we used nearest neighbor data imputation to complete missing aPTT observations. The feed forward ANN was trained using PyTorch's SGD method with Nestetorv momentum of 0.9, a batch size of 200, learning rate of 0.1, and 10 training epochs. The GRU and LSTM models were trained using PyTorch's Adam optimizer with a learning rate of 0.1, 50 training epochs, sequences that contain at most 48 transitions, and a batch size of 60 sequences.

\subsubsection{Comparison Results}

\begin{figure}
     \centering
     \begin{subfigure}[b]{0.49\textwidth}
         \centering
         \includegraphics[width=\textwidth]{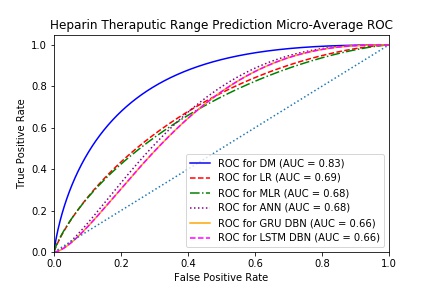}
         \caption{Micro average ROC}
         \label{fig:micro_roc}
     \end{subfigure}
     \hfill
     \begin{subfigure}[b]{0.49\textwidth}
         \centering
         \includegraphics[width=\textwidth]{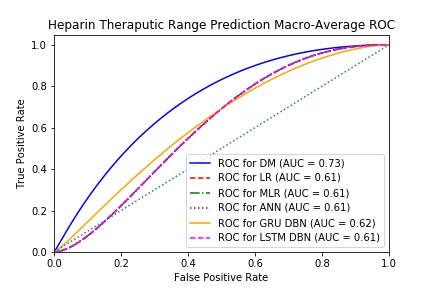}
         \caption{Macro average ROC}
         \label{fig:macro_roc}
     \end{subfigure}
        \caption{Micro and Macro average ROC comparison between different models for one vs all prediction of aPTT level. The x-axis measures the false positive rate and the y-axis measures the true positive rate of correctly predicting the aPTT level.}
        \label{ROC}
\end{figure}

Comparing the predictive performance of multi-class methods is less straight forward then in the case of binary classification since false positives and false negatives are harder to define. A common approach is to use a one vs all prediction evaluation, where the multi-class problem is essentially treated as being several binary classification problems each corresponding to one of the labels \citep{multiclass}. This approach allows us to create ROC curves for each model and get a sense for their predictive performance. One nuanced point however is that constructing ROC curves from one vs all prediction requires making assumptions on the true frequency of the classes. There are two common methods for this aggregation, the first is known as forming a micro-average that essentially assumes the true frequency of classes is the same as that observed in the data meaning that class imbalance will greatly influence the result \citep{multiclass}. The second approach is called a macro-average ROC, this method computes the metric independently for each class and then take the average across them, therefore it can be interpreted as equalizing the frequency of each class in the data \citep{multiclass}. To account for the advantages and disadvantages of each method we compute both the micro and macro-average smooth ROC curves for the different prediction methods and show the results in Figures \ref{fig:micro_roc} and  \ref{fig:macro_roc} respectively.
%

As can be seen in both plots, our  model outperforms all the other methods by about 0.1 AUC for both micro and macro-average ROC, and the patient heparin model ROC is strictly larger then all other ROC curves. This means that regardless of an imbalance in patient data, our model performs significantly better than other methods in terms of predicting future aPTT. All other methods seem to have similar performance with the RNN based methods seeming better in the case of the macro average ROC while the classical methods seem more effective in the micro average case. 

While the ROC gives a good comparative understanding of how the dynamic method ranks against other predictive methods, since it is constructed from one vs all metrics it does not provide the granularity for when our model is making prediction errors. To analyze this, we present the confusion matrix for our model in Table \ref{table:confusion_matrix}. 
Our dynamic model is effective in predicting when patients are sub-therapeutic and in therapeutic range while being less effective at predicting if patients are super-therapeutic. However, the case of patients being super therapeutic was rare in our data set which explains why the overall accuracy of the model is quite high. Moreover, when super-therapeutic values are miss-classified they tend to be classified as therapeutic. Since our model is performing a regression task, many of these values are predicted close to the upper bound of the therapeutic range. This means that in practice the predicted values can still be effective for dose setting and the results of the confusion matrix are an artifact of the conversion to classification.

\begin{table}
\begin{tabular}{p{1.4cm}|p{2.8cm}|p{2.8cm}|p{2.8cm}|p{2.8cm}|p{2.8cm}|}
\multicolumn{6}{c}{\qquad  \text{Predicted labels}} \\
\cline{2-6} & \textbf{Percent of Total} & \textbf{Sub-Therapeutic} & \textbf{Therapeutic} & \textbf{Super-Therapeutic} & \textbf{True Positive Rate }\\
\cline{2-6} &\textbf{sub-Therapeutic} & 37.54\% & 11.59\% & 0.69\% & \textbf{75.35\%} \\
\cline{2-6} Ground Truth & \textbf{Therapeutic} & 10.03\% & 26.12\% & 2.08\% & \textbf{68.33\%} \\
\cline{2-6} &\textbf{Super-Therapeutic} & 4.15\% & 6.23\% & 1.56\% & \textbf{13.04\%} \\
\cline{2-6} &\textbf{False Positive Rate} & \textbf{27.42\%} & \textbf{40.55\%} & \textbf{64.0\%} & \textbf{65.22\%} \\
\cline{2-6} 
\end{tabular}
\caption{Confusion matrix for dynamic model. Each entry represents the fraction of samples of a given ground truth label (row label) classified as another label (column label). The diagonal sum represents the model accuracy.}
\label{table:confusion_matrix}
\end{table}

\subsection{Adaptive dosing}
 In this section, we evaluate how the PTC-SGm algorithm introduced in Section \ref{sec:dose_opt} is able to maintain individual aPTT within a safe range in a simulated ICU environment, and compare its performance against a weight based dosing approach.

\subsubsection{Simulation setup}
We simulated a 240 hour long (10 day long) heparin treatment regimen for the 25 patients we obtained from the MIMIC III data set, with 10 replicates per patient and method. We assumed that each patient's individual heparin and aPTT dynamics followed the model we propose in Section \ref{sec:model_desc} with all unknown parameters estimated using joint MLE and the full records of the patients available in the data. We used 0 mean Laplace noise to simulate aPTT measurement noise, the variance of the noise term was estimated using the  sample variance of aPTT for each individual patient calculated from the data. 
%
%
%
%
%
 In a similar manner to real world heparin treatment, we assumed each of the dosing methods evaluated obtains a new noisy aPTT observation at 6 hour intervals. Then if the method has a learning component, that component would update its parameter estimates using this new observation (and all observations obtained prior in the simulation), followed by an optimization component (or other dose calculation component) that would then calculate the heparin dosing sequence for the next 6 hours. This process would then repeat at each 6 hour interval as new measurements are collected. In each simulation replicate, we assumed that the first 72 hours of treatment were identical to the treatment given in the MIMIC III data, and that after this time point the dosing method of interest would take over and administer doses to patients. Partially, this was to ensure that the sufficient richness condition was satisfied for algorithms with learning components, since this would ensure 12 aPTT measurements are available to each method initially. In application, a 3 day ramp up period could be shortened if additional lab measures are taken within a single day of the ICU stay. Though another interpretation of this is as part of a clinician in the loop deployment. Here, the clinical team sets the initial bolus when a patient enters the ICU, and then the automated dosing methods provide dosing recommendations after some initial lab tests are conducted. 
%
%
%
%
%
%
%

In total we evaluated 11 different adaptive dosing methods, of these 9 where variation on PTC-SGm with different implementation parameters. Using MLE across the full patient records, we found that the estimate of $\alpha$ was either $0.500$ or $0.707$ for most patients (a half life of 30 minutes or 2 hours respectively), therefore in this adaptive dosing experiment we specify that $\mathcal{A} = \{0.500, 0.707\}$ in the parameter estimation models of PTC-SGm. We further evenly discretized parameter space $\mathcal{B}$ and considered two setups where $|\mathcal{B}|=5$ and $|\mathcal{B}|=10$ so that we could evaluate both a 10 scenario and 20 scenario implementation of PTC-SGm that we refer to as PTC-SG10 and PTC-SG20 respectively. 
%
%
%
%
%
%
In addition to the two different versions of the PTC-SGm algorithm, we consider a similar dosing protocol where for each period we obtain the MLE of patient parameters with the measurements obtained to date and optimize doses assuming the MLE parameters to be true (this is equivalent to solving problem (\ref{eq:optimal_dosing})). We call this algorithm PTC-MLE. As discussed in Section \ref{dose_formulation}, the PTC algorithms involve using a loss function that can be represented by MILP constraints, and for all PTC-SG20, PTC-SG10, PTC-MLE we tested three different loss functions that will be discussed in the next section. We limited the maximum dose of heparin per hour that the optimization of the PTC algorithms can administer to be 3000 units to ensure patient safety. In addition we consider a weight-based protocol that simulates current clinical practice. For this method the dose of heparin administered to patients is determined both by their weight and by their risk of hemorrhaging. 
%
%
 For our simulation study, we follow the Ventura County Hospital Adult Heparin Drip Protocol (\cite{hirsh2008parenteral}). We note that the MIMIC III data set does not include each patient's risk of bleeding, so we use the estimated homeostasis aPTT $y_b$ from the full patient records to decide the risk of bleeding. If the patients has a low $y_b$ then they are at a lower risk of bleeding and a high $y_b$ means a higher risk of bleeding. We also consider a naive adaptive dosing policy. This policy uses our model and MLE to estimate individual patient parameters and to identify the therapeutic range for each patient. Instead of designing optimal doses by solving an optimization problem, the naive policy will increase the heparin dose for the next 6 hours by 200 units if the patient's aPTT level is sub-therapeutic, and decrease the dose by 200 units for the next 6 hours if they are super-therapeutic. If the patient is determined to be in the therapeutic range, the dose level from the past 6 hours is maintained into the next 6 hour period. 
The experiments were run on a laptop computer with 3.2GHz processor and 8GB RAM. To asses each dosing method, we measure the average time each method is able to maintain patient aPTT in the therapeutic range, the maximum deviation form that therapeutic range if the algorithm is unable to maintain it, and the running time for both predict and control phases. These metrics can be thought of as measuring efficacy, patient safety, and computational practicality respectively.

\subsubsection{Choice of loss function}
As discussed in Section \ref{sec:dose_opt}, our goal is to maintain patients' aPTT within a safe range and this is done through minimizing the expected posterior loss associated with a certain policy. In medical practice, the therapeutic aPTT range can be derived from the patient homeostasis aPTT $y_b$, namely $1.5y_b$ to $2.5y_b$ \citep{basu1972prospective}. In our dose planning phase during the PTC-SGm algorithm, we directly estimate the patient's homeostasis aPTT $y_b$ , and thus through optimization modeling techniques we can construct loss functions that penalize the objective when the patient's aPTT is not within the estimated safety range. Assumption \ref{loss_function} restricts the type of loss function we can use but it remains a mild assumption because in practice we can design many functions in this category. In this study, we consider the following loss functions:

\textbf{Indicator loss:} A straight forward loss function is the indicator function of whether the aPTT is staying within therapeutic range. The loss incurred at time $t$ is $\mathbb{1}_{\{y_t < 1.5y_b \, \text{or} \, y_t > 2.5y_b\}}$, meaning it penalizes the algorithm when the aPTT is not in therapeutic range by a constant loss of 1 and provides a loss of 0 no matter where in the therapeutic range the aPTT measure is. The indicator function can be described by MILP constraints.

\textbf{Absolute deviation from the therapeutic range:} Instead of having a constant penalty when the aPTT level is not in the therapeutic range, we can penalize it more as it deviates farther from the ideal range. This is a reasonable consideration because larger deviations imply larger levels of health risk. Accordingly, distance from current aPTT to the therapeutic band, i.e, $y_t - 2.5y_b$ if super-therapeutic, $1.5y_b- y_t$ if sub-therapeutic and 0 in safe range can also be used as the loss function. The loss incurred at time t can be expressed as $\mathbb{1}_{\{y_t < 1.5y_b \, or \, y_t > 2.5y_b\}} \cdot |y_t-2y_b|-\mathbb{1}_{\{y_t < 1.5y_b \, or \, y_t > 2.5y_b\}} \cdot 0.5y_b$. Again this loss can be represented with MILP constraints, and provides a linear penalty as aPTT deviates further from the therapeutic range.

\textbf{Absolute deviation from the median of therapeutic range:} Absolute deviation from the median therapeutic range provides a more aggressive penalty since it provides linear penalties for any deviation of aPTT from the median including those that are in the therapeutic range.  Although there is no evidence in medical literature proving that maintaining this median therapeutic aPTT is the best practice, this loss function forces the algorithm to closely match a specific aPTT level if one is proposed. The loss incurred at time $t$ is $|y_t - 2y_b|$. However, from a safety perspective, this loss can be quite advantageous since it penalizes the aPTT level from being close to the boundaries of the therapeutic range. Since we only have partial information of the aPTT and the patient's therapeutic range this can be seen as potentially more robust then the other losses that treat the whole therapeutic range equally. 
%

\begin{remark}
    In real world clinical practice, being sub-therapeutic and super-therapeutic could imply different levels of risk. To account for this with our model we can weigh deviation from upper therapeutic bound and lower therapeutic bound differently. For the purposes of our simulation however we weigh both equally.
\end{remark}

\subsubsection{Results and discussion} 
\begin{figure}
  \centering
  \captionsetup{justification=centering}
  \begin{subfigure}[t]{0.45\textwidth}
  \vskip 0pt
    \centering
    \includegraphics[width=\textwidth]{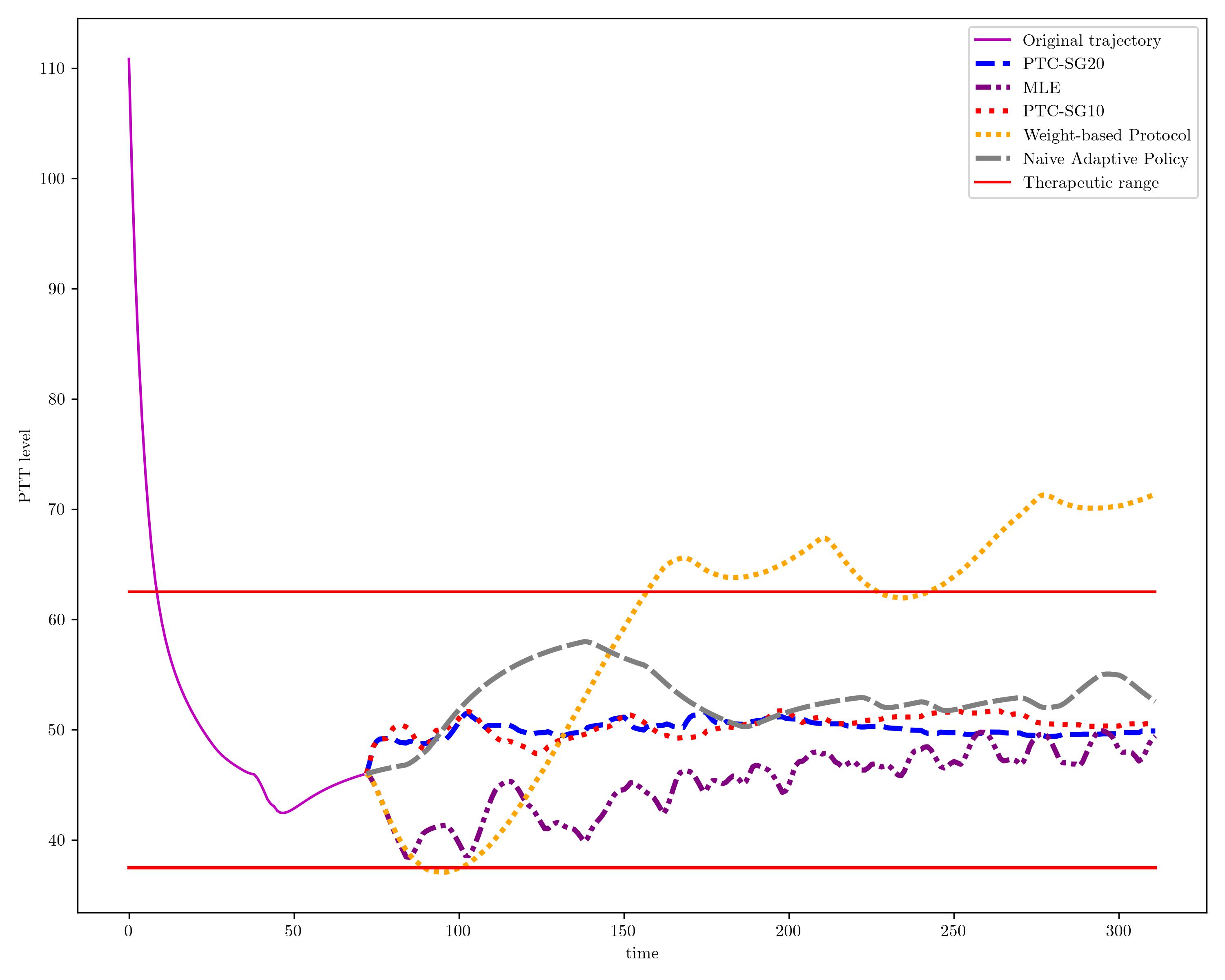}
    \caption{aPTT trajectories for an example patient under different dosing methods with the loss function of absolute deviation from the median of therapeutic range. The x-axis represents the time in hours in the ICU, and the y-axis represents the aPTT level. The simulated trajectories begin after time 72, the horizontal lines are the upper an lower bound of the therapeutic region.}
    \label{fig:aPTT_traj}
  \end{subfigure}
  \hfill
  \begin{subfigure}[t]{0.45\textwidth}
  \vskip 0pt
    \centering
    \includegraphics[width=\textwidth]{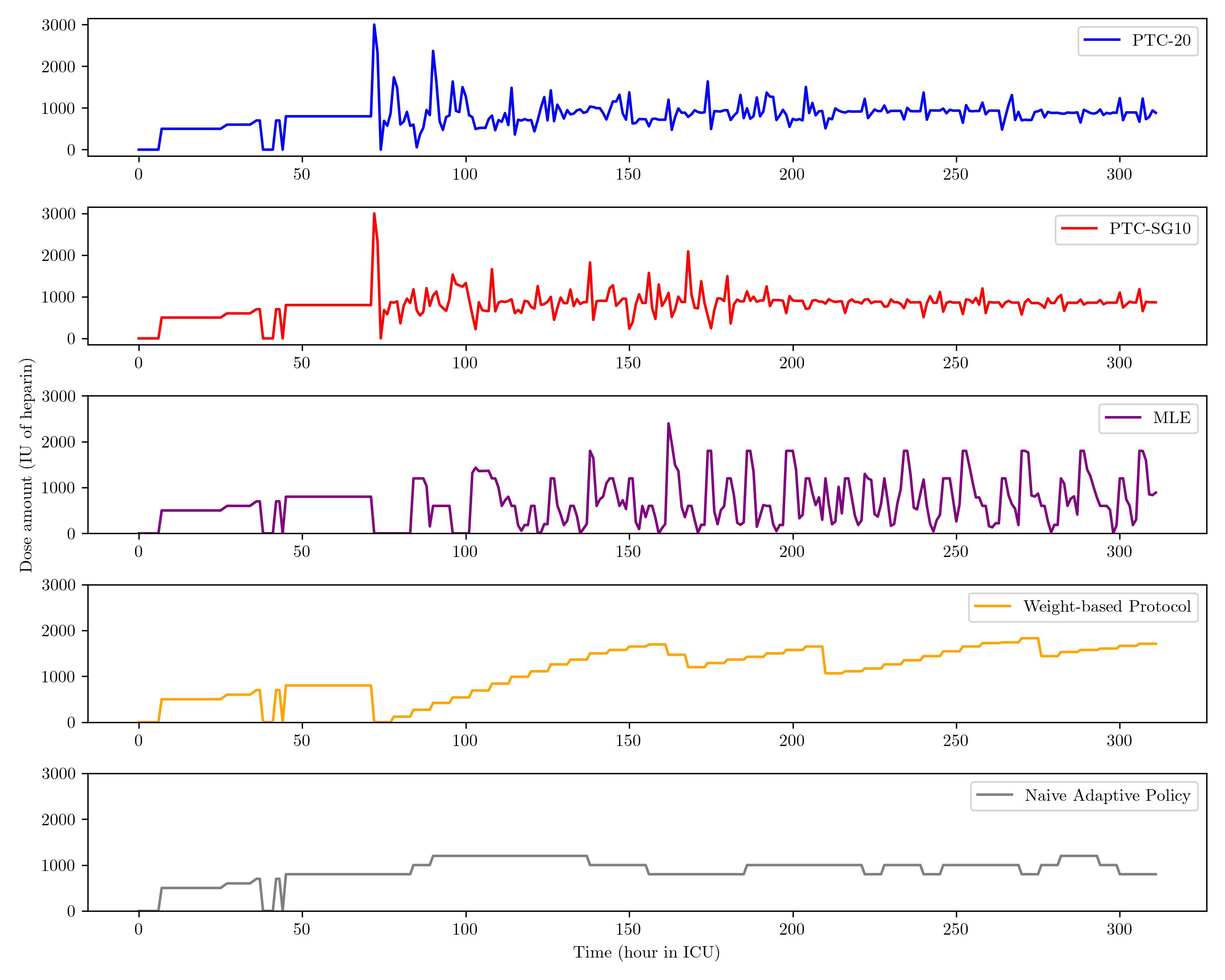}
    \caption{Dose amounts reccomended by each policy at each hour. The x-axis represents hours in the ICU the y-axis measures the units of heparin. From top to bottom the included plots are for the PTC-SG20, PTC-SG10, PTC-MLE, weight based dosing policy, and naive policy.}
    \label{fig:aPTT_dose}
  \end{subfigure}
    \caption{An example of personalized dosing policy including both aPTT trajectory and resulting heparin doses reccomended by each policy every hour.}
    \label{dosing}
\end{figure}

Figure \ref{fig:aPTT_traj} demonstrates an illustrative example of a single patient's aPTT trajectories under different dosing methods. Each of the PTC-SGm and PTC-MLE methods were implemented with the absolute deviation from the median of therapeutic range as the loss function. Figure \ref{fig:aPTT_dose} shows the corresponding dose sequences administered by each method. The x-axis of Figure \ref{fig:aPTT_traj} is time spent in ICU (in hours), and the y-axis is the patient's aPTT level at different times. The area between the two red lines is the ground truth therapeutic range where the patient is considered safe. According to the trajectory and the doses, when the patient was first admitted to the ICU, the aPTT is far above therapeutic, and the clinician had taken effective dosing steps to keep the aPTT under control. After 72 hours, the various methods being evaluated take charge of future doses. Both PTC-SG10 and PTC-SG20 algorithms behave similarly and keep the aPTT at median therapeutic level consistently. The PTC-MLE algorithm is not as stable as the scenario generation ones, producing more fluctuation in aPTT at the early stages. However, as the treatment continues, the fluctuation shrinks and PTC-MLE brings similar treatment effect as PTC-SGm. This indicates that a straight forward certainty equivalence technique only using the MLE estimate may not be as reliable in ensuring safety as the scenario generation based techniques. Moreover, while the weight-based protocol correctly increases the dose when the patient's aPTT is sub-therapeutic , it keeps dosing a sequence of large amount of heparin that put the patient in the dangerous super-therapeutic range. One reason for this is that the weight-based protocol cannot identify the ideal therapeutic level and targets at a higher aPTT than is apropriate for this patient. This is further reinforced as the naive dosing approach that uses our learning model to estimate the therapeutic range is capable of maintaining the patient in a safe range. However, much like the MLE approach, it provides less stable aPTT trajectories then the scenario generation approaches.

Figure \ref{fig:aPTT_dose} shows the initial doses given by the clinician and the amount of doses calculated by different dosing protocols at different times. Generally, the doses calculated by PTC-SG10 and PTC-SG20 overlap in both time and amount, indicating that for this patient the additional scenarios did not provide any benefit in treatment computation. PTC-MLE is not successful at the beginning because it fails to recognize that the patient's aPTT is going to decrease and  plans a sequences of 0 doses, which could be a safety concern. The naive approach however faces the opposite problem, since it can only modify the hourly dose slowly it provides steadily high doses that cause the aPTT to increase in the early stages, this again could be a safety concern leading to hemorrhaging. The weight based protocol however, does not adapt to the patient's data and constantly increases the heparin dose administered with periodic sudden drops. Overall, these dose plots suggest that both the predictive model that can estimate the therapeutic range, and the scenario generation method that considers model uncertainty could be key for patient safety.

%

For the patient in Figures \ref{fig:aPTT_traj} and \ref{fig:aPTT_dose}, all dosing methods other then the weight based protocol were able to keep the patient's aPTT within the therapeutic range for the entirety of the treatment.  While not all patients have a similarly perfect aPTT trajectory under the dosing protocol, the PTC-SGm algorithms achieve strong results accross all 25 patients. Table \ref{tb:dose_meth_comp} summarizes how the different methods perform on the whole patient population. We consider several metrics averaged over the 25 patients and their 40 dose planning cycles: time in control is the percentage of time the patients spent in therapeutic range during the 300-hour stay in ICU. Deviation from therapeutic range is the deviation of the patients' aPTT from their ideal range when they are not in the safe range. Predict time stands for the running time of the program to solve the scenario weight evaluation problem (Step 1 in Algorithm \ref{alg:scenario_generation}) or parameter estimation for PTC-MLE and the naive protocol. Control time stands for the running time of the program to solve the dose optimization problem. The method that achieves the most time in control is the PTC-SG10 algorithm  with the deviation form the median of the therapeutic range loss that keeps patients in therapeutic range for 87.7\% of time. Regarding different loss functions, the deviation from the median of therapeutic range has the best performance in terms of both time in control and deviation from the therapeutic range. This is because both the indicator of not staying in therapeutic range and absolute deviation from the therapeutic range loss functions will allow the patient's aPTT to stay close to the edge of therapeutic band (but within it) without penalty, which makes it more vulnerable to parameter estimation error and following doses calculated from very noisy observations. Regarding the number of scenarios used for PTC-SGm, it is surprising to see that PTC-SG10 sometimes outperforms PTC-SG20 in average time in control depending on the loss used. This is because when there are more scenarios, and the differences between the scenario with true parameters and another scenario are small, the algorithm will need more observations to adjust posterior probabilities for the different scenarios so that mass is centered on the scenario with the true parameters. Generally speaking, PTC-SGm converges slower when there are more scenarios. However, since the computation time for each planning cycle is within 1 minute, it seems that the running time should not be a concern for real time application. The predict time is much longer than the control time because the predict phase involves solving several profile likelihood optimization problems in parallel depending on the number of scenarios. 
%
%
Overall, the weight-based protocol is outperformed by most model-based methods in terms of time spent in the therapeutic range except for only two cases where a less robust loss functions and only the MLE is used without other scenarios. Interestingly, we note that even the naive adaptive policy is around 10\% better then the weight-based protocol in terms of time in the therapeutic range. This difference is mainly due to the naive policy using our predictive model to inform dosing decisions, thus showing that personalizing is key in providing treatment. However, PTC-SG10 and PTC-SG20 with the deviation from the median and deviation from the therapeutic range loss functions both significantly outperform the naive policy. This indicates that both the personalization of predictions to patients, as well as the design of the dose optimization method provides significant benefits to patient safety and health.

Through our simulation we found that the aPTT for some patients would not increase even though they were given large dose amounts. This is because these patients have an extremely low heparin coefficient $b$ meaning they are less responsive to heparin. Although heparin treatment is not effective for these patients, our model helps identify these patients and increase the chance that they are treated by other more effective means, unlike existing protocols that would simply increase their dose level without clinician input. 

\begin{figure}
    \centering
    \includegraphics[width=0.5\textwidth]{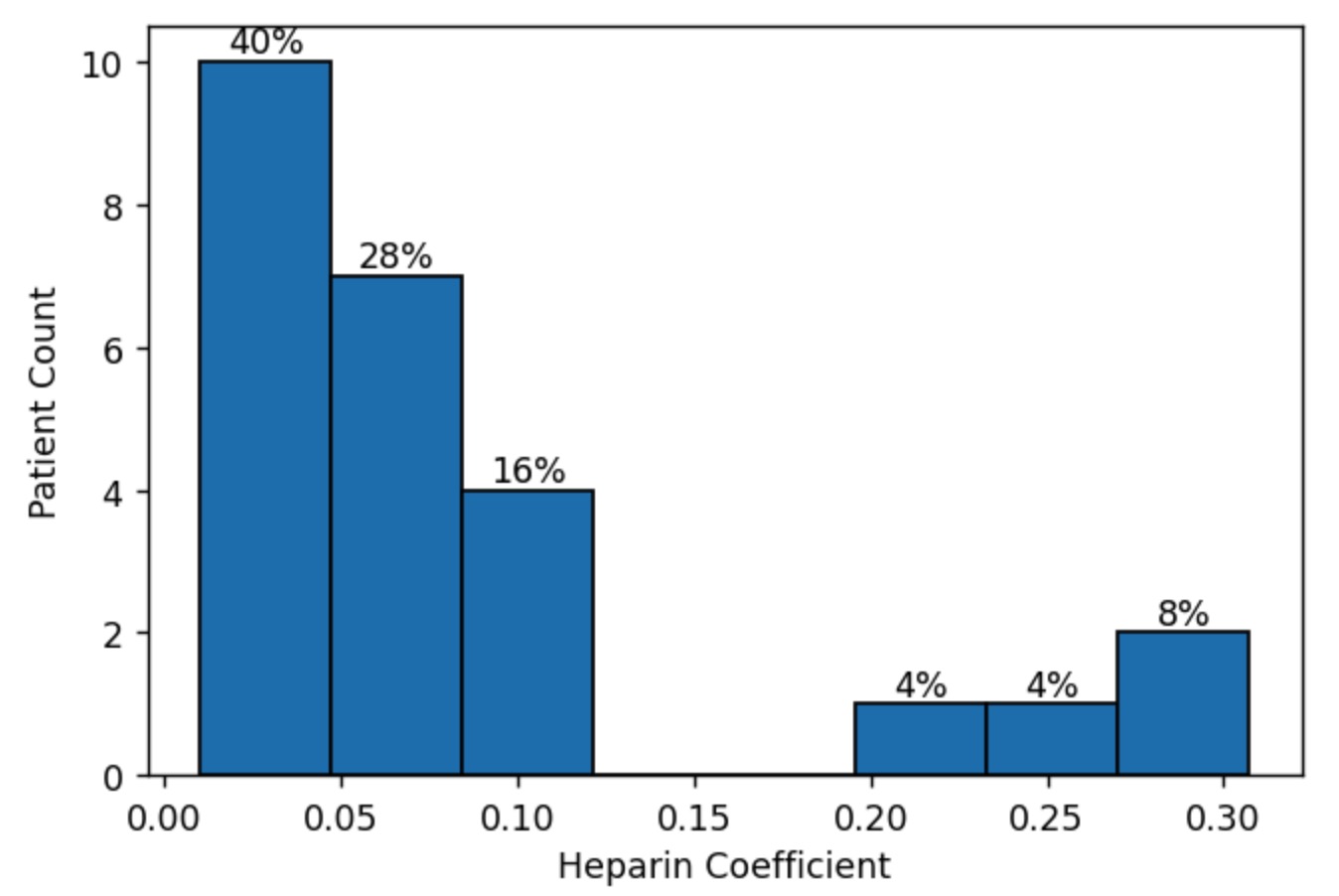}
    \caption{Heparin coefficient histogram computed from MLE estimates from individual patient data.}
    \label{fig:b_coefficient}
\end{figure}

\medskip
\begin{table}
\begin{tabular}{|p{6cm}|p{4cm}|p{1.6cm}|p{1.6cm}|p{1.3cm}|p{1.3cm}|}
\hline
\multicolumn{2}{|c|}{Dosing Protocol} &Time in Control & Deviation from therapeutic Range &Predict Time (s) &Control Time (s) \\
\hline
\multicolumn{2}{|c|}{Weight-based Protocol} & 55.6\% & 6.552 & 0 &0\\
\hline
\multicolumn{2}{|c|}{Naive Adaptive Dosing}  & 68.0\% & 2.639 & \textbf{2.507} & 0\\
\hline
\multirow{3}{*}{\shortstack[l]{PTC with Deviation \\from the Median of therapeutic \\Range as Loss Function}} & 20 Scenario Generation & 86.9\% & 0.447 & 12.750 & 0.168 \\
                     & 10 Scenario Generation & \textbf{87.7\%}  & \textbf{0.424} & 8.724 & 0.087 \\
                     & MLE & 69.9\% & 2.091 & 11.883 & 0.009\\
\hline
\multirow{3}{*}{\shortstack[l]{PTC with Deviation \\from therapeutic Range \\as Loss Function}} & 20 Scenario Generation & 82.4\% & 1.029 & 24.455 & 0.121 \\ 
                     & 10 Scenario Generation & 80.3\%  & 0.933 & 20.966 & 0.062 \\
                     & MLE & 53.0\% & 2.723 & 9.698 & 0.007\\
\hline
\multirow{3}{*}{\shortstack[l]{PTC with Indicator of \\staying in therapeutic Range \\or not as Loss Function}} & 20 Scenario Generation & 57.5\% & 5.238 & 7.926 & 0.040 \\
                     & 10 Scenario Generation &  59.6\% & 2.637 & 5.535 & 0.022 \\
                     & MLE & 36.2\% & 6.559 & 2.953 & \textbf{0.003}\\
\hline
\end{tabular}
\caption{Performance of dosing policies in the simulation experiment. Each entry represents an average across all 25 patients, at each of the 40 decision epochs, and 10 replicates. The best measure in each column is bolded.}
\label{tb:dose_meth_comp}
\end{table}

\subsection{Managerial Insights}
Our work provides insights not only for better  heparin dosing practices in ICUs, but also for the broader application of precision medicine in hospitals.
\begin{enumerate}
    \item \textit{Data driven methods can provide safer and more effective treatment then existing protocols.} Our experiment results show that even a naive adaptive approach that uses personalized estimates can outperform the current weight-based dosing protocol and that the best PTC-SGm algorithm outperforms the weight-based dosing protocol by about 32.1\% in terms of keeping patients in the therapeutic range. Also, all the algorithms can be run within in a minute. This implies that data-driven  decision-making methods may have better accuracy then existing protocols and can be deployed at scale. Clinicians could use the recommendations and predictions from these models to enhance their outcomes when facing complex decisions a cross a larger number of patients.
    \item \textit{Personalization is key to patient safety especially in the case of precision medicine.} Our estimates for patient parameters show that there is significant patient heterogeneity in heparin sensitivity as shown by Figure \ref{fig:b_coefficient}. This individual difference cannot be identified by a general machine learning method that learns population level parameters from a group of patients. Moreover, estimating the patient specific therapeutic ranges is key to ensuring each individual is kept safe from adverse effects such as hemorrhaging or clotting, unlike the existing protocols that target a single range. Therefore, in choosing or developing models for medical decision making, personalizing is essential for effective treatment planning and patient safety.

    \item \textit{Optimization methods developed for precision medicine must account for likely adverse outcomes and not just mean outcome to ensure patient safety.} Our PTC-SGm algorithm leverages different combinations of patient parameters and adapts to patient's true parameters gradually as more observations are obtained so it is robust to large observation errors. In contrast to this, methods that are not adaptive and rely too much on initial data with small size are vulnerable to observation errors as shown by our computational results. Therefore, in developing methods that aim to implement precision medicine  patient safety must be emphasized. Methods should not only be accurate in expectation, but also need to be robust to possible errors during implementation, and must be adaptive to changing health conditions.
\end{enumerate}
\section{Conclusion}
In this paper, we developed a novel model-based framework for personalized heparin dosing in ICUs. Overcoming the issue of uncertainties in both the patient's states and dynamics faced by traditional model-free approaches, our pharmacokinetic model makes accurate prediction of patients' states. We then proposed a scenario generation algorithm called PTC-SGm that optimizes future doses using this model. Every component in the framework can be computed using commercial MILP solvers. We provided a theoretical analysys that showed that PTC-SGm provides asymptotically optimal dosing sequences. We then concluded with computational experiments that validate the performance of our methods.


\newpage
\begin{APPENDICES}
\section{Proofs of Propositions in Text}
\label{app:proofs}
\proof{Proof of Proposition \ref{prop:z_bx}} Using the substitution in the premise of the proof, consider the one step dynamics function of $x_{t+1}=h\left(x_t\right)+u_t$. Since $b>0$, we can multiply through by $b$ and see that $z_{t+1}=b h\left(x_t\right)+b u_t$. Using the structure of $h$ as a piecewise linear function we can rewrite $b h\left(x_t\right)$ by multiplying through each piece by $b$ as:
\begin{equation}
b h\left(x_t\right)= \begin{cases}\alpha z_t, & \text { for } x_t \leq \frac{k}{1-\alpha} \\ z_t-b k & \text { for } x_t>\frac{k}{1-\alpha}\end{cases}.
\end{equation}
To adjust the break point to be in terms of $z_t$, we can multiply both sides of the inequality and show that the break points occur at $z_t \leq \frac{b k}{1-\alpha}$ and $z_t \geq \frac{b k}{1-\alpha}$. Defining $c=b k$, completes the reformulation. The reformulation of the terms relating to $y_{t+1}$ follows by direct substitution.\halmos

\proof{Proof of Proposition \ref{prop:alpha_z_reform}} Let $\bar{z}_t$ and $\bar{\alpha}=\alpha_{\bar{i}} \in \mathcal{A}$ be a solution to the original problem. Consider the set of vectors $z'_t=\bar{z}_t$, $\iota'_{i}=1$ for $i=\bar{i}$, $\iota'_{i}=0$ for $i \neq\bar{i}$, $w'_{it} = \bar{\alpha}z_t$ for $i=\bar{i}$, $w'_{it} = 0$ for $i \neq \bar{i}$ and $w'_t = \sum_{i = 1}^m w'_{it}$. It is easy to check that the set of vectors $z'_t,\iota'_i,w'_{it},w'_t$ is a feasible solution to the reformulated problem. On the other hand, let $z'_t,\iota'_i,w'_{it},w'_t$ be a set of feasible solution to the reformulated problem where $\iota'_{\bar{i}}=1$. Consider the set of vectors $\bar{\alpha} = \alpha_{\bar{i}} \in \mathcal{A}$, $\bar{z}_t=bu_{t-1}+\alpha_{\bar{i}}\bar{z}_{t-1}$ if $\bar{z}_{t-1} \leq \frac{c}{1-\alpha}$, $\bar{z}_t=bu_{t-1}+\bar{z}_{t-1-c}$ if $\bar{z}_{t-1} > \frac{c}{1-\alpha}$ for $t=1,2,...,T$, $\bar{z}_0 = 0$. It is easy to check that the set of vectors $\bar{\alpha}, \bar{z}_t$ is a set of feasible solution to the original problem. We conclude the proposition. \halmos
\endproof

\proof{Proof of Proposition \ref{prop:MILP_solver}} Let $\bar{z}_t, \bar{w}_t$ be a set of feasible solution to the original problem. Consider the set of vectors $z'_t = \bar{z}_t$, $w'_t=\bar{w}_t$, $\nu'_t = 1$ if $\bar{z}_t \geq c + \bar{w}_t$, $\nu'_t = 0$ if $\bar{z}_t \leq c + \bar{w}_t$. It is easy to check that the set of vectors $z'_t, w'_t, \nu't_t$ is a feasible solution to the reformulated problem. On the other hand, let $z'_t, w'_t, \nu't_t$ be a feasible solution to the reformulated problem. Consider the set of vectors $\bar{z}_t = z'_t$, $\bar{w}_t = w'_t$, if $\nu'_t = 0$, then $\bar{z}_t \leq c + \bar{w}_t$ and $\bar{z}_{t+1} = \bar{w}_t+bu_t$; if $\nu'_t = 1$, then $\bar{z}_t \geq c + \bar{w}_t$ and $\bar{z}_{t+1} = \bar{z}_t-c$. We verify that the set of vectors $\bar{z}_t, \bar{w}_t$ is a feasible solution to the original problem. \halmos
\endproof

\proof{Proof of Proposition \ref{prop:lp_sub_prob}} Observe that for fixed $ \bar{\alpha} \in \mathcal{A}, \bar{k} \in \mathcal{K},\bar{b} \in \mathcal{B}$ the entire sequence of $z_t$ is determined (since by assumption $z_0 = 0$). For simplicity of notation let us call $\bar{z}_t = z_t(\bar{\alpha},\bar{k}, \bar{b}, \{u_\tau\}_{\tau=0}^t)$, i.e. the value at time $t$ of the modified state trajectory given these fixed starting conditions  and let us call $\overrightarrow{z}(\bar{\alpha},\bar{k},\bar{b})= [0, \bar{z}_1,...,\bar{z}_T ]^\top$ . Hence we can rewrite the optimization problem in \eqref{eq:sub_prob} as follows:
\begin{equation}
    \begin{aligned}
   \mathcal{L}(\overrightarrow{z}(\bar{a},\bar{k},\bar{b})) = \max_{\{y_t,y_{bt}\}_{t=0}^T,y_b} & \sum_{t \in n(T)} p(\tilde{y}_t | y_t) \\
\text{subject to: } &\\
& y_{t+1} = \gamma_1 (y_t - y_{bt}) + y_{bt} + \bar{z}_{t+1}, & \forall t \in \{0,...,T-1\}, \\
&y_{b,t+1} = \gamma_2  y_{bt} + \gamma_3y_{t} + \gamma_4 y_b, & \forall t \in \{0,...,T-1\},\\ 
&  y_t,y_{bt},y_b \in \mathcal{Y}, & \forall t \in \{0,...,T\}.
    \end{aligned}
\end{equation}

Observe that we have now rewritten $\mathcal{L}$ as the value function of a convex problem with respect to the right hand side of its constraints, this proves the desired results. \halmos
\endproof

\proof{Proof of Proposition \ref{prop:master}}
First note that Problem \eqref{eq:constraind_mle} can be expressed as: $\max_{a \in \mathcal{A}, k\in \mathcal{K},b\in\mathcal{B}} \mathcal{L}(\overrightarrow{z}(a,k,b)) $ by simply substituting the appropriate value. We can alternatively express this problem in the following form using the definition of the maximum:
\begin{equation}
\label{eq:master_form}
\begin{aligned}
    \max_{a , k ,b,\ell}& \ \ell \\
    &\text{subject to: } \\
    & \ell \leq \mathcal{L}(\overrightarrow{z}(a,k,b) ),\\
    & a \in \mathcal{A}, k \in \mathcal{K},b\in \mathcal{B}.
    \end{aligned}
\end{equation}
Next, by Proposition \ref{prop:lp_sub_prob} we know that $\mathcal{L}$ is the value function of a convex optimization problem for fixed $\overrightarrow{z}$, therefore using strong duality we can consider the the Lagrangian dual of \eqref{eq:sub_prob}.

\begin{equation}
\label{eq:lagrange}
    \begin{aligned}
    \mathcal{L}(\overrightarrow{z}(\bar{\alpha},\bar{k},\bar{b}))= \min_{\lambda} \overrightarrow{z}(\bar{\alpha},\bar{k},\bar{b})^\top \overrightarrow{\lambda} + g(\{\lambda_t\}_{t=0}^{n-1}).
    \end{aligned}
\end{equation}

Note that $|g(\{\lambda_t\}_{t=0}^{n-1})| \leq \infty $ for any fixed sequence $\{\lambda_t\}_{t=0}^{n-1}$, this is because the optimization problem in \eqref{eq:dual_sub} will always have $y_{b,t} = 0$ for all $t = \{0,...,n\}$ as a feasible solution and by Assumption 2 $p(\cdot)$ is log concave. Therefore Problem \eqref{eq:lagrange} will always be feasible for fixed $\overrightarrow{z}$. Therefore we need to consider two alternatives, either \eqref{eq:lagrange} obtains an optimal solution, or it is unbounded meaning that Subproblem \eqref{eq:sub_prob} is infeasible (and hence the given sequence $\overrightarrow{z}$ is infeasible for the original problem).  In the case where a feasible solution is obtained we can directly substitute the Lagrangian objective into  \eqref{eq:master_form}
giving us the first set of constraints from \eqref{eq:master}. 

Therefore all that is left is to address the case of infeasibility. If Problem \eqref{eq:sub_prob} is infeasible then so is the corresponding feasibility problem:
\begin{equation}
\label{eq:feasibility}
    \begin{aligned}
        \min_{\{y_t,y_{bt}\}_{t=0}^T,y_b} \ & 0 \\
        & \text{subject to:} \\
        & y_{t+1} = \gamma_1 (y_t - y_{bt}) + y_{bt} + z_{t+1}, & \forall t \in \{0,...,T-1\}, \\
&y_{b,t+1} = \gamma_2  y_{bt} + \gamma_3y_{t} + \gamma_4 y_b, & \forall t \in \{0,...,T-1\},\\
    & y_{bt},y_{t},y_{b} \in \mathcal{Y}, & \forall t \in \{0,...,T\}.
    \end{aligned}
\end{equation}
Observe that this feasibility problem is a linear program. Consider the constraint set \begin{equation}
    \begin{aligned}
        &y_{b,t+1} = \gamma_2  y_{bt} + \gamma_3y_{t} + \gamma_4 y_b, & \forall t \in \{0,...,T-1\},\\
    & y_{bt},y_{t},y_{b} \in \mathcal{Y}, & \forall t \in \{0,...,T\}.
    \end{aligned}
\end{equation}
Note that these constraints are always feasible, by Assumption \ref{as:alpha_set_assump} we know interval $\mathcal{Y}$ is nonempty, so we can choose some point $\bar{y}\in\mathcal{Y}$ and set $y_b = \bar{y}$, $y_{b0} = \bar{y}$ and $y_t = \bar{y}, \forall t$. Since each following $y_{bt}$ is a convex combination of points in $\mathcal{Y}$, and $\mathcal{Y}$ is a convex set by assumption these starting conditions constitute a feasible point. Hence to ensure the fixed $\bar{\alpha}, \bar{k}, \bar{b}$ (and the resulting $\bar{z}$ are feasible), we just need to add the first constraint in \eqref{eq:feasibility} to \eqref{eq:master_form}. Let $S_1 = \{\overrightarrow{z}: y_{t+1} = \gamma_1 (y_t - y_{bt}) + y_{bt} + \bar{z}_{t+1}, y_{b,t+1} = \gamma_2  y_{bt} + \gamma_3y_{t} + \gamma_4 y_b, \forall t \in \{0,...,n-1\} \textit{ for some } \theta_y \in \Theta_y$\} , $S_2 = \{\overrightarrow{z}: 0 = \overrightarrow{z}(\bar{\alpha},\bar{k},\bar{b})^\top \overrightarrow{\mu} + r_y(\{\mu_t\}_{t=0}^{n-1}) \text{ for all }  \overrightarrow{\mu} \in [-1,1]^T \}$, we need to show that $S_1=S_2$.

Let $\overrightarrow{z}$ be an arbitrary point in $S_1$. $\bar{z}_t \in S_2$ since the first equality defining $S_1$ can also be written as $\overrightarrow{z}^{\top} \overrightarrow{\mu} + \sum_{t=0}^{n-1} \mu_t(\gamma_1 (y_t - y_{bt}) + y_{bt} - y_{t+1})$ for all $\overrightarrow{\mu} \in [-1,1]^T $. The converse can be demonstrated with the help of duality theory. Suppose $\overrightarrow{z} \in S_2$, then we have $\smash{\displaystyle\min_{\mu \in \mathcal{R}^n}\{ \overrightarrow{z}(\bar{\alpha},\bar{k},\bar{b})^\top \overrightarrow{\mu} + r_y(\{\mu_t\}_{t=0}^{n-1})\}=0}$ since the scaling of $\overrightarrow{\mu}$ does not influence the equality. Consider the feasibility problem
\begin{equation}
\label{eq:feasi_dual}
\begin{aligned}
    \max_{\{y_t,y_{bt}\}_{t=0}^T,y_b}& \ 0 \\
    &\text{s.t. } \\
    & y_{t+1} = \gamma_1 (y_t - y_{bt}) + y_{bt} + \bar{z}_{t+1}, \forall t \in \{0,...,T-1\},\\
    & y_{b,t+1} = \gamma_2  y_{bt} + \gamma_3y_{t} + \gamma_4 y_b, & \forall t \in \{0,...,T-1\}, \\
    & y_{bt},y_{t},y_{b} \in \mathcal{Y}, & \forall t \in \{0,...,T\}.
    \end{aligned}
\end{equation}
 $\smash{\displaystyle\min_{\mu \in \mathcal{R}^n}\{ \overrightarrow{z}(\bar{\alpha},\bar{k},\bar{b})^\top \overrightarrow{\mu} + r_y(\{\mu_t\}_{t=0}^{n-1})\}=0}$ asserts that the dual of this feasibility problem with respect to the first constraint has optimal value 0. Therefore, this optimization problem must be feasible, i.e, $\overrightarrow{z} \in S_1$.
\halmos 

\endproof

\proof{Proof of Proposition \ref{prop:benders_convergence}}
By our assumption, the sets $\mathcal{B} \times \mathcal{K}$ and $\Theta_y$ are compact. The objective function $\sum_{t \in n(T)} \log p(\tilde{y}_t| y_t)$ is concave and continuous on $\mathcal{B} \times \mathcal{K}$ and $\Theta_y$. The subproblem is always feasible as shown in Proposition \ref{prop:master}. Therefore by Theorem 2.5 in \cite{geoffrion1972generalized}, for each fixed $\alpha$ and $\forall \epsilon > 0$, the decomposition procedure terminates finitely. \halmos
\endproof

\proof{Proof of Corollary \ref{corol:bayesian_consistency}} The proof follows by a similar method presented in \cite{mintz2017behavioral}. Let $(\theta^*_x, \theta^*_y)$ be the patient's true parameters, and observe that

\begin{equation}
    \begin{aligned}
        &\log  \hat{p}\big(\theta_x,\theta_y  | \{\tilde{y}_t\}_{t\in n(T)}, \{u_t\}_{t=0}^T \big)\\
&=
    \log \big( \hat{p}(\theta^*_x,\theta^*_y | \{\tilde{y}_t\}_{t\in n(T)},\{u_t\}_{t=0}^T) \cdot \frac{\hat{p}( \theta_x,\theta_y| \{\tilde{y}_t\}_{t\in n(T)} ,\{u_t\}_{t=0}^T)}{\hat{p}(\theta^*_x,\theta^*_y |\{\tilde{y}_t\}_{t\in n(T)},\{u_t\}_{t=0}^T )} \big)\\
&=
    \sum_{t\in n(T)} \log \frac{p_{\epsilon}\left(\tilde{y}_{t}- \bar{y}_{t}\right)}{p_{\epsilon}\left(\tilde{y}_{t}- y^*_{t}\right)} +
    \sum_{t=1}^T \log\frac{ p(\bar{y}_{t-1},\bar{y}_{bt}, \bar{x}_t, \bar{\theta}_y|\bar{y}_t)}{ p(y^*_{t-1},y^*_{bt}, x^*_t, \theta^*_y|y^*_t)} + 
    \sum_{t=1}^T \log \frac{p(\bar{x}_{t-1}, u_t, \bar{\theta}_x|\bar{x}_t )}{p(x^*_{t-1}, u_t, \theta^*_x|x^*_t)} \\
&+ \log \hat{p}(\theta^*_x,\theta^*_y | \{\tilde{y}_t\}_{t\in n(T)},\{u_t\}_{t=0}^T)
    \end{aligned}
\end{equation}

where $x^*_t$, $y^*_t$ are the states under true dynamics and the patient's initial conditions $(\theta^*_x, \theta^*_y)$, and $\bar{x}_t$, $\bar{y}_t$ are the states under estimated parameters $(\theta_x,\theta_y)$. Notice that the dynamics $y_t = f(x_t,y_{t-1},\theta_y)$ and $x_t  = g(x_{t-1},u_t,\theta_x)$ hold for both true states and estimated states, so $p(y_{t-1},y_{bt},x_t,\theta_y|y_t )$ and $p(x_{t-1},u_t,\theta_x|x_t)$ are degenerate. Therefore $\log\frac{ p(\bar{y}_{t-1},\bar{y}_{bt},\bar{x}_t, \bar{\theta}_y|\bar{y}_t)}{ p(y^*_{t-1},y^*_{bt}, x^*_t, \theta^*_y|y^*_t)} = \log \frac{p(\bar{x}_{t-1}, u_t, \bar{\theta}_x|\bar{x}_t )}{p(x^*_{t-1}, u_t, \theta^*_x|x^*_t )} = \log1 = 0$. Also, $\log \hat{p}(\theta^*_x,\theta^*_y | \{\tilde{y}_t\}_{t\in n(T)},\{u_t\}_{t=0}^T) \leq 0$ since $ \hat{p}(\theta^*_x,\theta^*_y | \{\tilde{y}_t\}_{t\in n(T)},\{u_t\}_{t=0}^T) \in [0,1]$ Then by Assumption 1 we have $\max _{\mathcal{E}(\delta)} \log p(\{\tilde{y}_t\}_{t\in n(T)} |\bar{\theta}_x,\bar{\theta}_y) \rightarrow -\infty$ for any $\delta > 0$ almost surely.  While this shows the point-wise convergence of the 
estimate, we can derive its uniform convergence. For any $\delta > 0$ as $T \rightarrow \infty $ we have that

\begin{equation}
    \label{uniform_convergence}
    \begin{gathered}
\hat{p}(\mathcal{E}(\delta))|\{\tilde{y}_t\}_{t\in n(T)}, \{u_t\}_{t=0}^T)=\int_{\mathcal{E}(\delta)} \hat{p}(\theta_x,\theta_y)|\{\tilde{y}_t\}_{t\in n(T)}, \{u_t\}_{t=0}^T) \times d \theta_x \times d \theta_y \leq \\
\operatorname{volume}(\Theta_x \times \Theta_y) \cdot \max _{\mathcal{E}(\delta)} \hat{p}\left(\theta_x,\theta_y)|\{\tilde{y}_t\}_{t\in n(T)}, \{u_t\}_{t=0}^T\right) \rightarrow 0,
    \end{gathered}
\end{equation}
which completes the proof. \halmos
\endproof

\proof{Proof of Corollary \ref{corol:estimator_consistency}}
Consider the events $E_1=\left\{(\hat{\theta}_x, \hat{\theta}_y, \notin \mathcal{B}\left(\theta^*_x, \theta^*_y, \delta\right)\right\}$ and $E_2=\left\{\max _{\mathcal{E}(\delta)} \hat{p}(\theta_x,\theta_y \mid\tilde{y}_t\}_{t\in n(T)}, \{u_t\}_{t=0}^T) \geq \max _{(\theta_x, \theta_y) \in \mathcal{B}(x_T^*, \theta_T^*, \delta)} \hat{p}(x_T, \theta_T \mid\tilde{y}_t\}_{t\in n(T)}, \{u_t\}_{t=0}^T\right\}$ where $\mathcal{E}(\delta)$ is defined as before for some $\delta>0$. Then observe that $E_1 \subset E_2$, therefore $\mathbb{P}\left(E_1\right) \leq \mathbb{P}\left(E_2\right)$. By Proposition 2 as $T \rightarrow \infty, \mathbb{P}\left(E_2\right) \rightarrow 0$ hence $\mathbb{P}\left(E_1\right) \rightarrow 0$. Thus the result of the corollary follows. \halmos
\endproof

\proof{Proof of Corollary \ref{coro:dose}}
Apply the reformulation introduced in Section \ref{sec:milp_form_mle}  (\ref{eq:optimal_dosing}). Notice that $k,y_b,y_0,y_{b0},\{u_t\}_{t=0}^{T+n}$ belong to an affine term in mixed integer linear constraints. Standard results \citep{ralphs2014value} imply that $\varphi_{\bar{\alpha},\bar{b}}$, as the value function of a MILP, is lower semicontinuous with respect to $k,y_b,y_0,y_{b0},\{u_t\}_{t=0}^{T+n}$. \halmos
\endproof

\proof{Proof of Proposition \ref{algorithm_asymptotic1}} Corollary (\ref{corol:bayesian_consistency}) implies that 
$\frac{\omega_{\alpha^*,b^*}}{\sum\limits_{\bar{\alpha},\bar{b} \in \bar{\mathcal{A}} \times \bar{\mathcal{B}}} \omega_{\bar{\alpha},\bar{b}}} \xrightarrow{p} 1$, and for $\forall (\alpha',b') \neq (\alpha^*,b^*)$ $\frac{\omega_{\alpha',b'}}{\sum\limits_{\bar{\alpha},\bar{b} \in \bar{\mathcal{A}} \times \bar{\mathcal{B}}} \omega_{\bar{\alpha},\bar{b}}} \xrightarrow{p}  0 $. By property of profile likelihood estimation \citep{murphy2000profile} we have that $\tau(\alpha^*,b^*) \xrightarrow{p} (k^*,y_b^*,y_0^*,y_{b0}^*)$. According to continuous mapping theorem \citep{van2000asymptotic}, we have that $\varphi_{\alpha^*,b^*}(\tau(\alpha^*,b^*),\{u_t\}_{t=0}^{T+n}) \xrightarrow{p} \varphi_{\alpha^*,b^*}\left(k^*,y_b^*,y_0^*,y_{b0}^*,\{u_t\}_{t=0}^{T+n}\right)$. Putting all above together, we get the desired result. \halmos
\endproof

\proof{Proof of Proposition \ref{algorithm_asymptotic2}}
Definition 4.2 from \cite{vogel2003continuous} suggests that convergence in probability implies lower semicontinuous approximation in probability. Therefore 
with Proposition \ref{algorithm_asymptotic1}, for fixed dose sequence $\{u_t\}_{t=T+1}^{T+n}$, we get a point-wise result $\frac{\sum\limits_{\bar{\alpha},\bar{b} \in \mathcal{A} \times \mathcal{B}} \varphi_{\bar{\alpha},\bar{b}}(\tau(\bar{\alpha},\bar{b}),\{u_t\}_{t=0}^{T+n})\omega_{\bar{\alpha},\bar{b}}}{\sum\limits_{\bar{\alpha},\bar{b} \in \mathcal{A} \times \mathcal{B}} \omega_{\bar{\alpha},\bar{b}}} \underset{\{u_t\}_{t=T+1}^{T+n}}{\stackrel{l-p r o b}{\longrightarrow}} \varphi_{\alpha^*,b^*}\left(k^*,y_b^*,y_0^*,y_{b0}^*,\{u_t\}_{t=0}^{T+n}\right))$. This is valid for all $\{u_t\}_{t=T+1}^{T+n} \in \mathcal{U}^n$. Based on the relationship between point-wise lowersemicontinuous and uniform lowersemicontinous stated in Proposition 5.1 from \cite{vogel2003continuous} we conclude the uniform lowersemicontinuous result. \halmos
\endproof

\proof{Proof of Theorem \ref{asymptotic_doses}}
The result follows by applying Proposition \ref{algorithm_asymptotic2} combined with Theorem 4.3 from \citep{vogel2003continuous}. \halmos
\endproof
\end{APPENDICES}
\newpage
\newcommand{\newblock}{}
\bibliography{reference}

\begin{thebibliography}{}

\bibitem[\protect\astroncite{Amiral et~al.}{2021}]{heparin_1}
Amiral, J., Amiral, C., and Dunois, C. (2021).
\newblock Optimization of heparin monitoring with anti-fxa assays and the
  impact of dextran sulfate for measuring all drug activity.
\newblock {\em Biomedicines}, 9(6):700.

\bibitem[\protect\astroncite{{\AA}str{\"o}m and
  Wittenmark}{2013}]{aastrom2013adaptive}
{\AA}str{\"o}m, K.~J. and Wittenmark, B. (2013).
\newblock {\em Adaptive control}.
\newblock Courier Corporation.

\bibitem[\protect\astroncite{Aswani et~al.}{2019}]{aswani2019behavioral}
Aswani, A., Kaminsky, P., Mintz, Y., Flowers, E., and Fukuoka, Y. (2019).
\newblock Behavioral modeling in weight loss interventions.
\newblock {\em European journal of operational research}, 272(3):1058--1072.

\bibitem[\protect\astroncite{Ayer et~al.}{2012}]{pdosing_7}
Ayer, T., Alagoz, O., and Stout, N.~K. (2012).
\newblock Or forum—a pomdp approach to personalize mammography screening
  decisions.
\newblock {\em Operations Research}, 60(5):1019--1034.

\bibitem[\protect\astroncite{Basu et~al.}{1972}]{basu1972prospective}
Basu, D., Gallus, A., Hirsh, J., and Cade, J. (1972).
\newblock A prospective study of the value of monitoring heparin treatment with
  the activated partial thromboplastin time.
\newblock {\em New England Journal of Medicine}, 287(7):324--327.

\bibitem[\protect\astroncite{Beal}{1983}]{MM3}
Beal, S.~L. (1983).
\newblock Computation of the explicit solution to the michaelis-menten
  equation.
\newblock {\em Journal of pharmacokinetics and biopharmaceutics},
  11(6):641--657.

\bibitem[\protect\astroncite{Belotti et~al.}{2016}]{bigM}
Belotti, P., Bonami, P., Fischetti, M., Lodi, A., Monaci, M.,
  Nogales-G{\'o}mez, A., and Salvagnin, D. (2016).
\newblock On handling indicator constraints in mixed integer programming.
\newblock {\em Computational Optimization and Applications}, 65(3):545--566.

\bibitem[\protect\astroncite{Bertsekas}{2012}]{POMDP_1}
Bertsekas, D. (2012).
\newblock {\em Dynamic programming and optimal control: Volume I}, volume~1.
\newblock Athena scientific.

\bibitem[\protect\astroncite{Bertsimas
  et~al.}{2020}]{bertsimas2020personalized}
Bertsimas, D., Orfanoudaki, A., and Weiner, R.~B. (2020).
\newblock Personalized treatment for coronary artery disease patients: a
  machine learning approach.
\newblock {\em Health care management science}, 23:482--506.

\bibitem[\protect\astroncite{Bickel and Doksum}{2015}]{bickel2015mathematical}
Bickel, P.~J. and Doksum, K.~A. (2015).
\newblock {\em Mathematical statistics: basic ideas and selected topics,
  volumes I-II package}.
\newblock Chapman and Hall/CRC.

\bibitem[\protect\astroncite{Bjornsson and Nash}{1986}]{aPTT_linear}
Bjornsson, T.~D. and Nash, P.~V. (1986).
\newblock Variability in heparin sensitivity of aptt reagents.
\newblock {\em American journal of clinical pathology}, 86(2):199--204.

\bibitem[\protect\astroncite{Bonifonte et~al.}{2022}]{bonifonte2022analytics}
Bonifonte, A., Ayer, T., and Haaland, B. (2022).
\newblock An analytics approach to guide randomized controlled trials in
  hypertension management.
\newblock {\em Management Science}, 68(9):6634--6647.

\bibitem[\protect\astroncite{Boyd et~al.}{2004}]{boyd2004convex}
Boyd, S., Boyd, S.~P., and Vandenberghe, L. (2004).
\newblock {\em Convex optimization}.
\newblock Cambridge university press.

\bibitem[\protect\astroncite{Brunet et~al.}{2008}]{brunet2008pharmacodynamics}
Brunet, P., Simon, N., Opris, A., Faure, V., Lorec-Penet, A.-M., Portugal, H.,
  Dussol, B., and Berland, Y. (2008).
\newblock Pharmacodynamics of unfractionated heparin during and after a
  hemodialysis session.
\newblock {\em American journal of kidney diseases}, 51(5):789--795.

\bibitem[\protect\astroncite{Bull et~al.}{1975}]{heparin_problem}
Bull, B.~S., Korpman, R.~A., Huse, W.~M., and Briggs, B.~D. (1975).
\newblock Heparin therapy during extracorporeal circulation: I. problems
  inherent in existing heparin protocols.
\newblock {\em The Journal of thoracic and cardiovascular surgery},
  69(5):674--684.

\bibitem[\protect\astroncite{{\c{C}}elik et~al.}{2015}]{POMDP_3}
{\c{C}}elik, M., Ergun, {\"O}., and Keskinocak, P. (2015).
\newblock The post-disaster debris clearance problem under incomplete
  information.
\newblock {\em Operations Research}, 63(1):65--85.

\bibitem[\protect\astroncite{Cho et~al.}{2014}]{cho2014learning}
Cho, K., Van~Merri{\"e}nboer, B., Gulcehre, C., Bahdanau, D., Bougares, F.,
  Schwenk, H., and Bengio, Y. (2014).
\newblock Learning phrase representations using rnn encoder-decoder for
  statistical machine translation.
\newblock {\em arXiv preprint arXiv:1406.1078}.

\bibitem[\protect\astroncite{Conforti et~al.}{2014}]{conforti2014integer}
Conforti, M., Cornu{\'e}jols, G., Zambelli, G., et~al. (2014).
\newblock {\em Integer programming}, volume 271.
\newblock Springer.

\bibitem[\protect\astroncite{Cook}{2010}]{cook2010anticoagulation}
Cook, B.~W. (2010).
\newblock Anticoagulation management.
\newblock In {\em Seminars in interventional radiology}, volume~27, pages
  360--367. {\copyright} Thieme Medical Publishers.

\bibitem[\protect\astroncite{Cornish-Bowden}{2015}]{MM1}
Cornish-Bowden, A. (2015).
\newblock One hundred years of michaelis--menten kinetics.
\newblock {\em Perspectives in Science}, 4:3--9.

\bibitem[\protect\astroncite{Coroian and Hauser}{2015}]{pdosing_6}
Coroian, D.~C. and Hauser, K. (2015).
\newblock Learning stroke treatment progression models for an mdp clinical
  decision support system.
\newblock In {\em Proceedings of the 2015 SIAM International Conference on Data
  Mining}, pages 676--684. SIAM.

\bibitem[\protect\astroncite{Craig et~al.}{1987}]{craig1987adaptive}
Craig, J.~J., Hsu, P., and Sastry, S.~S. (1987).
\newblock Adaptive control of mechanical manipulators.
\newblock {\em The International Journal of Robotics Research}, 6(2):16--28.

\bibitem[\protect\astroncite{Davenport}{2011}]{davenport2011optimization}
Davenport, A. (2011).
\newblock Optimization of heparin anticoagulation for hemodialysis.
\newblock {\em Hemodialysis International}, 15:S43--S48.

\bibitem[\protect\astroncite{Dogan et~al.}{2021}]{dogan2021regret}
Dogan, I., Shen, Z.-J.~M., and Aswani, A. (2021).
\newblock Regret analysis of learning-based mpc with partially-unknown cost
  function.
\newblock {\em arXiv preprint arXiv:2108.02307}.

\bibitem[\protect\astroncite{Eagle}{1984}]{POMDP_2}
Eagle, J.~N. (1984).
\newblock The optimal search for a moving target when the search path is
  constrained.
\newblock {\em Operations research}, 32(5):1107--1115.

\bibitem[\protect\astroncite{Eikelboom and Hirsh}{2006}]{heparin_aPTT}
Eikelboom, J.~W. and Hirsh, J. (2006).
\newblock Monitoring unfractionated heparin with the aptt: time for a fresh
  look.
\newblock {\em Thrombosis and haemostasis}, 96(11):547--552.

\bibitem[\protect\astroncite{Embretson and Reise}{2013}]{embretson2013item}
Embretson, S.~E. and Reise, S.~P. (2013).
\newblock {\em Item response theory}.
\newblock Psychology Press.

\bibitem[\protect\astroncite{Garcia and Fern{\'a}ndez}{2012}]{saferl_5}
Garcia, J. and Fern{\'a}ndez, F. (2012).
\newblock Safe exploration of state and action spaces in reinforcement
  learning.
\newblock {\em Journal of Artificial Intelligence Research}, 45:515--564.

\bibitem[\protect\astroncite{Gaskett}{2003}]{saferl_2}
Gaskett, C. (2003).
\newblock Reinforcement learning under circumstances beyond its control.

\bibitem[\protect\astroncite{Geibel and Wysotzki}{2005}]{saferl_3}
Geibel, P. and Wysotzki, F. (2005).
\newblock Risk-sensitive reinforcement learning applied to control under
  constraints.
\newblock {\em Journal of Artificial Intelligence Research}, 24:81--108.

\bibitem[\protect\astroncite{Geoffrion}{1972}]{geoffrion1972generalized}
Geoffrion, A.~M. (1972).
\newblock Generalized benders decomposition.
\newblock {\em Journal of optimization theory and applications},
  10(4):237--260.

\bibitem[\protect\astroncite{Ghassemi et~al.}{2014}]{multilogit}
Ghassemi, M.~M., Richter, S.~E., Eche, I.~M., Chen, T.~W., Danziger, J., and
  Celi, L.~A. (2014).
\newblock A data-driven approach to optimized medication dosing: a focus on
  heparin.
\newblock {\em Intensive care medicine}, 40(9):1332--1339.

\bibitem[\protect\astroncite{Greene}{2003}]{greene2003econometric}
Greene, W.~H. (2003).
\newblock {\em Econometric analysis}.
\newblock Pearson Education India.

\bibitem[\protect\astroncite{Guillet et~al.}{2003}]{heparin2}
Guillet, B., Simon, N., Sampol, J.~J., Lorec-Penet, A.-M., Portugal, H.,
  Berland, Y., Dussol, B., and Brunet, P. (2003).
\newblock Pharmacokinetics of the low molecular weight heparin enoxaparin
  during 48 h after bolus administration as an anticoagulant in haemodialysis.
\newblock {\em Nephrology Dialysis Transplantation}, 18(11):2348--2353.

\bibitem[\protect\astroncite{Hagan et~al.}{1997}]{hagan1997neural}
Hagan, M.~T., Demuth, H.~B., and Beale, M. (1997).
\newblock {\em Neural network design}.
\newblock PWS Publishing Co.

\bibitem[\protect\astroncite{Hajjar and Alagoz}{2023}]{hajjar2023personalized}
Hajjar, A. and Alagoz, O. (2023).
\newblock personalized disease screening decisions considering a chronic
  condition.
\newblock {\em Management Science}, 69(1):260--282.

\bibitem[\protect\astroncite{Hauskrecht and Fraser}{2000}]{pdosing_4}
Hauskrecht, M. and Fraser, H. (2000).
\newblock Planning treatment of ischemic heart disease with partially
  observable markov decision processes.
\newblock {\em Artificial intelligence in medicine}, 18(3):221--244.

\bibitem[\protect\astroncite{Hirsh et~al.}{2001}]{hirsh2001guide}
Hirsh, J., Anand, S.~S., Halperin, J.~L., and Fuster, V. (2001).
\newblock Guide to anticoagulant therapy: Heparin: a statement for healthcare
  professionals from the american heart association.
\newblock {\em Circulation}, 103(24):2994--3018.

\bibitem[\protect\astroncite{Hirsh et~al.}{2008}]{hirsh2008parenteral}
Hirsh, J., Bauer, K.~A., Donati, M.~B., Gould, M., Samama, M.~M., and Weitz,
  J.~I. (2008).
\newblock Parenteral anticoagulants: American college of chest physicians
  evidence-based clinical practice guidelines.
\newblock {\em Chest}, 133(6):141S--159S.

\bibitem[\protect\astroncite{Hochreiter and
  Schmidhuber}{1997}]{hochreiter1997long}
Hochreiter, S. and Schmidhuber, J. (1997).
\newblock Long short-term memory.
\newblock {\em Neural computation}, 9(8):1735--1780.

\bibitem[\protect\astroncite{Hou and Jin}{2013}]{model_free}
Hou, Z. and Jin, S. (2013).
\newblock {\em Model free adaptive control}.
\newblock CRC press Boca Raton, FL.

\bibitem[\protect\astroncite{Jambhekar and Breen}{2009}]{jambhekar2009basic}
Jambhekar, S.~S. and Breen, P.~J. (2009).
\newblock {\em Basic pharmacokinetics}, volume~76.
\newblock Pharmaceutical press London.

\bibitem[\protect\astroncite{Janczak and Grishin}{2006}]{unknown_state}
Janczak, D. and Grishin, Y. (2006).
\newblock State estimation of linear dynamic system with unknown input and
  uncertain observation using dynamic programming.
\newblock {\em Control and Cybernetics}, 35(4):851--862.

\bibitem[\protect\astroncite{Jiang and Powell}{2015}]{jiang2015approximate}
Jiang, D.~R. and Powell, W.~B. (2015).
\newblock An approximate dynamic programming algorithm for monotone value
  functions.
\newblock {\em Operations research}, 63(6):1489--1511.

\bibitem[\protect\astroncite{Johnson et~al.}{2016}]{mimic3}
Johnson, A.~E., Pollard, T.~J., Shen, L., Lehman, L.-w.~H., Feng, M., Ghassemi,
  M., Moody, B., Szolovits, P., Anthony~Celi, L., and Mark, R.~G. (2016).
\newblock Mimic-iii, a freely accessible critical care database.
\newblock {\em Scientific data}, 3(1):1--9.

\bibitem[\protect\astroncite{Karakala and Tolwani}{2016}]{weight2}
Karakala, N. and Tolwani, A. (2016).
\newblock We use heparin as the anticoagulant for crrt.
\newblock In {\em Seminars in dialysis}, volume~29, pages 272--274. Wiley
  Online Library.

\bibitem[\protect\astroncite{Kaut and Stein}{2003}]{scenario}
Kaut, M. and Stein, W. (2003).
\newblock {\em Evaluation of scenario-generation methods for stochastic
  programming}.
\newblock Humboldt-Universit{\"a}t zu Berlin,
  Mathematisch-Naturwissenschaftliche Fakult{\"a}t~….

\bibitem[\protect\astroncite{Keskinocak and Savva}{2020}]{keskinocak2020review}
Keskinocak, P. and Savva, N. (2020).
\newblock A review of the healthcare-management (modeling) literature published
  in manufacturing \& service operations management.
\newblock {\em Manufacturing \& Service Operations Management}, 22(1):59--72.

\bibitem[\protect\astroncite{Kong et~al.}{2017}]{deeparch}
Kong, N., Liu, X., Liu, C., Lian, J., and Wang, H. (2017).
\newblock Deep architecture for heparin dosage prediction during continuous
  renal replacement therapy.
\newblock In {\em 2017 36th Chinese Control Conference (CCC)}, pages
  11166--11171. IEEE.

\bibitem[\protect\astroncite{Koppejan and Whiteson}{2011}]{saferl_4}
Koppejan, R. and Whiteson, S. (2011).
\newblock Neuroevolutionary reinforcement learning for generalized control of
  simulated helicopters.
\newblock {\em Evolutionary intelligence}, 4(4):219--241.

\bibitem[\protect\astroncite{Landefeld
  et~al.}{1987}]{landefeld1987identification}
Landefeld, C.~S., Cook, E.~F., Flatley, M., Weisberg, M., and Goldman, L.
  (1987).
\newblock Identification and preliminary validation of predictors of major
  bleeding in hospitalized patients starting anticoagulant therapy.
\newblock {\em The American journal of medicine}, 82(4):703--713.

\bibitem[\protect\astroncite{Lee et~al.}{2015}]{lee2015applying}
Lee, E., Lavieri, M.~S., Volk, M.~L., and Xu, Y. (2015).
\newblock Applying reinforcement learning techniques to detect hepatocellular
  carcinoma under limited screening capacity.
\newblock {\em Health care management science}, 18:363--375.

\bibitem[\protect\astroncite{Lee et~al.}{2018}]{lee2018outcome}
Lee, E.~K., Wei, X., Baker-Witt, F., Wright, M.~D., and Quarshie, A. (2018).
\newblock Outcome-driven personalized treatment design for managing diabetes.
\newblock {\em Interfaces}, 48(5):422--435.

\bibitem[\protect\astroncite{Li}{2017}]{li2017deep}
Li, Y. (2017).
\newblock Deep reinforcement learning: An overview.
\newblock {\em arXiv preprint arXiv:1701.07274}.

\bibitem[\protect\astroncite{Li et~al.}{2016}]{meta2}
Li, Y., Monine, M., Huang, Y., Swann, P., Nestorov, I., and Lyubarskaya, Y.
  (2016).
\newblock Quantitation and pharmacokinetic modeling of therapeutic antibody
  quality attributes in human studies.
\newblock In {\em MAbs}, volume~8, pages 1079--1087. Taylor \& Francis.

\bibitem[\protect\astroncite{Liu et~al.}{2017}]{pdosing_2}
Liu, Y., Logan, B., Liu, N., Xu, Z., Tang, J., and Wang, Y. (2017).
\newblock Deep reinforcement learning for dynamic treatment regimes on medical
  registry data.
\newblock In {\em 2017 IEEE international conference on healthcare informatics
  (ICHI)}, pages 380--385. IEEE.

\bibitem[\protect\astroncite{Maier et~al.}{2021}]{pdosing_1}
Maier, C., Hartung, N., Kloft, C., Huisinga, W., and de~Wiljes, J. (2021).
\newblock Reinforcement learning and bayesian data assimilation for
  model-informed precision dosing in oncology.
\newblock {\em CPT: pharmacometrics \& systems pharmacology}, 10(3):241--254.

\bibitem[\protect\astroncite{Mart{\'\i}n~H and Lope}{2009}]{saferl_6}
Mart{\'\i}n~H, J.~A. and Lope, J.~d. (2009).
\newblock Learning autonomous helicopter flight with evolutionary reinforcement
  learning.
\newblock In {\em International Conference on Computer Aided Systems Theory},
  pages 75--82. Springer.

\bibitem[\protect\astroncite{McAvoy}{1979}]{MM2}
McAvoy, T. (1979).
\newblock Pharmacokinetic modeling of heparin and its clinical implications.
\newblock {\em Journal of Pharmacokinetics and Biopharmaceutics},
  7(4):331--354.

\bibitem[\protect\astroncite{Mintz et~al.}{2017}]{mintz2017behavioral}
Mintz, Y., Aswani, A., Kaminsky, P., Flowers, E., and Fukuoka, Y. (2017).
\newblock Behavioral analytics for myopic agents.
\newblock {\em arXiv preprint arXiv:1702.05496}.

\bibitem[\protect\astroncite{Murphy and Van~der
  Vaart}{2000}]{murphy2000profile}
Murphy, S.~A. and Van~der Vaart, A.~W. (2000).
\newblock On profile likelihood.
\newblock {\em Journal of the American Statistical Association},
  95(450):449--465.

\bibitem[\protect\astroncite{Narasimhan et~al.}{2016}]{multiclass}
Narasimhan, H., Pan, W., Kar, P., Protopapas, P., and Ramaswamy, H.~G. (2016).
\newblock Optimizing the multiclass f-measure via biconcave programming.
\newblock In {\em 2016 IEEE 16th international conference on data mining
  (ICDM)}, pages 1101--1106. IEEE.

\bibitem[\protect\astroncite{Nemati et~al.}{2016}]{Nemati2016OptimalMD}
Nemati, S., Ghassemi, M.~M., and Clifford, G.~D. (2016).
\newblock Optimal medication dosing from suboptimal clinical examples: A deep
  reinforcement learning approach.
\newblock {\em 2016 38th Annual International Conference of the IEEE
  Engineering in Medicine and Biology Society (EMBC)}, pages 2978--2981.

\bibitem[\protect\astroncite{Oudemans-van Straaten et~al.}{2011}]{weight1}
Oudemans-van Straaten, H.~M., Kellum, J.~A., and Bellomo, R. (2011).
\newblock Clinical review: Anticoagulation for continuous renal replacement
  therapy-heparin or citrate?
\newblock {\em Critical Care}, 15(1):1--9.

\bibitem[\protect\astroncite{Paszke et~al.}{2017}]{paszke2017automatic}
Paszke, A., Gross, S., Chintala, S., Chanan, G., Yang, E., DeVito, Z., Lin, Z.,
  Desmaison, A., Antiga, L., and Lerer, A. (2017).
\newblock Automatic differentiation in pytorch.

\bibitem[\protect\astroncite{Pedregosa et~al.}{2011}]{scikit-learn}
Pedregosa, F., Varoquaux, G., Gramfort, A., Michel, V., Thirion, B., Grisel,
  O., Blondel, M., Prettenhofer, P., Weiss, R., Dubourg, V., Vanderplas, J.,
  Passos, A., Cournapeau, D., Brucher, M., Perrot, M., and Duchesnay, E.
  (2011).
\newblock Scikit-learn: Machine learning in {P}ython.
\newblock {\em Journal of Machine Learning Research}, 12:2825--2830.

\bibitem[\protect\astroncite{Qin et~al.}{2014}]{unknown_dynamics}
Qin, C., Zhang, H., and Luo, Y. (2014).
\newblock Online optimal tracking control of continuous-time linear systems
  with unknown dynamics by using adaptive dynamic programming.
\newblock {\em International Journal of Control}, 87(5):1000--1009.

\bibitem[\protect\astroncite{Ralphs and Hassanzadeh}{2014}]{ralphs2014value}
Ralphs, T.~K. and Hassanzadeh, A. (2014).
\newblock On the value function of a mixed integer linear optimization problem
  and an algorithm for its construction.
\newblock {\em COR@ L Technical Report 14T--004}.

\bibitem[\protect\astroncite{Roberts}{2003}]{meta1}
Roberts, S.~A. (2003).
\newblock Drug metabolism and pharmacokinetics in drug discovery.
\newblock {\em Current opinion in drug discovery \& development}, 6(1):66--80.

\bibitem[\protect\astroncite{Rockafellar and
  Wets}{2009}]{rockafellar2009variational}
Rockafellar, R.~T. and Wets, R. J.-B. (2009).
\newblock {\em Variational analysis}, volume 317.
\newblock Springer Science \& Business Media.

\bibitem[\protect\astroncite{Sato et~al.}{2001}]{saferl_1}
Sato, M., Kimura, H., and Kobayashi, S. (2001).
\newblock Td algorithm for the variance of return and mean-variance
  reinforcement learning.
\newblock {\em Transactions of the Japanese Society for Artificial
  Intelligence}, 16(3):353--362.

\bibitem[\protect\astroncite{Secretariat}{2009}]{secretariat2009point}
Secretariat, M.~A. (2009).
\newblock Point-of-care international normalized ratio (inr) monitoring devices
  for patients on long-term oral anticoagulation therapy: an evidence-based
  analysis.
\newblock {\em Ontario health technology assessment series}, 9(12):1.

\bibitem[\protect\astroncite{Shani et~al.}{2013}]{POMDP_4}
Shani, G., Pineau, J., and Kaplow, R. (2013).
\newblock A survey of point-based pomdp solvers.
\newblock {\em Autonomous Agents and Multi-Agent Systems}, 27(1):1--51.

\bibitem[\protect\astroncite{Shen and Winkelmayer}{2012}]{shen2012use}
Shen, J.~I. and Winkelmayer, W.~C. (2012).
\newblock Use and safety of unfractionated heparin for anticoagulation during
  maintenance hemodialysis.
\newblock {\em American journal of kidney diseases}, 60(3):473--486.

\bibitem[\protect\astroncite{Shi et~al.}{2021}]{shi2021timing}
Shi, P., Helm, J.~E., Deglise-Hawkinson, J., and Pan, J. (2021).
\newblock Timing it right: Balancing inpatient congestion vs. readmission risk
  at discharge.
\newblock {\em Operations Research}, 69(6):1842--1865.

\bibitem[\protect\astroncite{Skandari and Shechter}{2021}]{skandari2021patient}
Skandari, M.~R. and Shechter, S.~M. (2021).
\newblock Patient-type bayes-adaptive treatment plans.
\newblock {\em Operations Research}, 69(2):574--598.

\bibitem[\protect\astroncite{Sui et~al.}{2015}]{safebandit1}
Sui, Y., Gotovos, A., Burdick, J., and Krause, A. (2015).
\newblock Safe exploration for optimization with gaussian processes.
\newblock In {\em International conference on machine learning}, pages
  997--1005. PMLR.

\bibitem[\protect\astroncite{Sun et~al.}{2017}]{safebandit2}
Sun, W., Dey, D., and Kapoor, A. (2017).
\newblock Safety-aware algorithms for adversarial contextual bandit.
\newblock In {\em International Conference on Machine Learning}, pages
  3280--3288. PMLR.

\bibitem[\protect\astroncite{Sutton and Barto}{2018}]{sutton2018reinforcement}
Sutton, R.~S. and Barto, A.~G. (2018).
\newblock {\em Reinforcement learning: An introduction}.
\newblock MIT press.

\bibitem[\protect\astroncite{Van~der Vaart}{2000}]{van2000asymptotic}
Van~der Vaart, A.~W. (2000).
\newblock {\em Asymptotic statistics}, volume~3.
\newblock Cambridge university press.

\bibitem[\protect\astroncite{Vogel and Lachout}{2003}]{vogel2003continuous}
Vogel, S. and Lachout, P. (2003).
\newblock On continuous convergence and epi-convergence of random functions.
  part i: Theory and relations.
\newblock {\em Kybernetika}, 39(1):75--98.

\bibitem[\protect\astroncite{Vozikis and Goulionis}{2009}]{pdosing_5}
Vozikis, A. and Goulionis, J.~E. (2009).
\newblock Medical decision making for patients with parkinson disease under
  average cost criterion.
\newblock {\em Australia and New Zealand health policy}, 6(1).

\bibitem[\protect\astroncite{Wang}{2022}]{wang2022optimal}
Wang, J. (2022).
\newblock Optimal sequential multiclass diagnosis.
\newblock {\em Operations Research}, 70(1):201--222.

\bibitem[\protect\astroncite{Wang et~al.}{2022}]{wang2022reliable}
Wang, J., Gao, R., and Zha, H. (2022).
\newblock Reliable off-policy evaluation for reinforcement learning.
\newblock {\em Operations Research}.

\bibitem[\protect\astroncite{Wilhelmsson and
  Lins}{1984}]{wilhelmsson1984heparin}
Wilhelmsson, S. and Lins, L. (1984).
\newblock Heparin elimination and hemostasis in hemodialysis.
\newblock {\em Clinical nephrology}, 22(6):303--306.

\bibitem[\protect\astroncite{Wolsey and Nemhauser}{1999}]{wolsey1999integer}
Wolsey, L.~A. and Nemhauser, G.~L. (1999).
\newblock {\em Integer and combinatorial optimization}, volume~55.
\newblock John Wiley \& Sons.

\bibitem[\protect\astroncite{Wu and Suen}{2022}]{wu2022optimizing}
Wu, C.-C. and Suen, S.-c. (2022).
\newblock Optimizing diabetes screening frequencies for at-risk groups.
\newblock {\em Health Care Management Science}, 25(1):1--23.

\bibitem[\protect\astroncite{Yu et~al.}{2021a}]{yu2021reinforcement}
Yu, C., Liu, J., Nemati, S., and Yin, G. (2021a).
\newblock Reinforcement learning in healthcare: A survey.
\newblock {\em ACM Computing Surveys (CSUR)}, 55(1):1--36.

\bibitem[\protect\astroncite{Yu et~al.}{2021b}]{pdosing_3}
Yu, C., Liu, J., Nemati, S., and Yin, G. (2021b).
\newblock Reinforcement learning in healthcare: A survey.
\newblock {\em ACM Computing Surveys (CSUR)}, 55(1):1--36.

\bibitem[\protect\astroncite{Zhang et~al.}{2012}]{zhang2012optimization}
Zhang, J., Denton, B.~T., Balasubramanian, H., Shah, N.~D., and Inman, B.~A.
  (2012).
\newblock Optimization of prostate biopsy referral decisions.
\newblock {\em Manufacturing \& Service Operations Management}, 14(4):529--547.

\end{thebibliography}




\end{document}